\documentclass[11pt,namelimits,sumlimits,a4paper]{amsart}
\usepackage{cite}
\usepackage{comment}

\usepackage{amssymb,amsmath}
\usepackage[mathscr]{eucal}
\usepackage{slashed}

\usepackage[latin1]{inputenc}
\usepackage[T1]{fontenc}
\usepackage[UKenglish]{babel}
\usepackage{amsfonts}
\usepackage{fancyhdr}
\usepackage{amsmath}
\usepackage{amsthm}
\usepackage{amsmath,amscd}
\usepackage{latexsym}
\usepackage{cite}
\usepackage{amssymb,amsmath}

\usepackage[mathscr]{eucal}
\usepackage{slashed}
\usepackage{times,psfrag}
\usepackage{mathtools}
\usepackage{hyperref}
\usepackage{multirow}
\usepackage{longtable}
\usepackage{bm}
\usepackage{fancyhdr}
\usepackage{float}

\usepackage{esint}

\usepackage{layout}

\usepackage{lmodern}
\usepackage[rgb]{xcolor}
\usepackage[draft,author={Lorenzo Foscolo}]{pdfcomment}

\usepackage{enumitem}

\numberwithin{equation}{section}
\newtheorem{theorem}[equation]{Theorem}
\newtheorem*{theorem*}{Theorem}
\newtheorem{lemma}[equation]{Lemma}
\newtheorem{prop}[equation]{Proposition}

\theoremstyle{definition}
\newtheorem{definition}[equation]{Definition}

\theoremstyle{remark}
\newtheorem{remark}[equation]{Remark}
\newtheorem*{remark*}{Remark}

\newcommand{\ie}{\emph{i.e.} }

\newcommand{\cf}{\emph{cf.} }

\newcommand{\beq}{\begin{equation}}
\newcommand{\eeq}{\end{equation}}
\newcommand{\bea}{\begin{eqnarray}}
\newcommand{\eea}{\end{eqnarray}}

\newcommand{\C}{\mathbb{C}}

\newcommand{\R}{\mathbb{R}}

\newcommand{\Z}{\mathbb{Z}}

\newcommand{\T}{\mathbb{T}}
\newcommand{\HH}{\mathbb{H}}

\newcommand{\PP}{\mathbb{P}}
\newcommand{\Sph}{\mathbb{S}}
\newcommand{\ra}{\rightarrow}

\newcommand{\vol}{\operatorname{Vol}}
\newcommand{\dvol}{\operatorname{dv}}

\newcommand{\triple}[1]{\boldsymbol{\underline{#1}}}

\def\co{\colon\thinspace}

\begin{document}

\title{ALF gravitational instantons and collapsing Ricci-flat metrics on the $K3$ surface}

\author[L.~Foscolo]{Lorenzo~Foscolo}
\address{Mathematics Department, State University of New York at Stony Brook}
\email{lorenzo.foscolo@stonybrook.edu}

\maketitle

\begin{abstract}
We construct large families of new collapsing hyperk\"ahler metrics on the K3 surface. The limit space is a flat Riemannian $3$--orbifold $T^3/\Z_2$. Away from finitely many exceptional points the collapse occurs with bounded curvature. There are at most $24$ exceptional points where the curvature concentrates, which always contains the $8$ fixed points of the involution on $T^3$. The geometry around these points is modelled by ALF gravitational instantons: of dihedral type ($D_k$) for the fixed points of the involution on $T^3$ and of cyclic type ($A_k$) otherwise.

The collapsing metrics are constructed by deforming approximate hyperk\"ahler metrics obtained by gluing ALF gravitational instantons to a background (incomplete) $S^1$--invariant hyperk\"ahler metric arising from the Gibbons--Hawking ansatz over a punctured $3$--torus.

As an immediate application to submanifold geometry, we exhibit hyperk\"ahler metrics on the K3 surface that admit a strictly stable minimal sphere which cannot be holomorphic with respect to any complex structure compatible with the metric.
\end{abstract}

\section{Introduction}

Soon after Yau's proof of the Calabi Conjecture showed that the (smooth $4$--manifold underlying a complex) K3 surface carries K\"ahler Ricci-flat metrics, Gibbons and Pope \cite{Gibbons:Pope} suggested the construction (further explored by Page in \cite{Page:K3}) of explicit approximately Ricci-flat metrics on Kummer surfaces. They considered the quotient of a flat $4$--torus by an involution and resolved the $16$ orbifold singularities by gluing in copies of the Eguchi--Hanson metric. This \emph{Kummer construction} was later made rigorous by Topiwala \cite{Topiwala} and LeBrun--Singer \cite{LeBrun:Singer:Kummer} using twistor methods. Much more recently Donaldson \cite{Donaldson:Kummer} gave a different proof using analysis, closest to the approach taken in the current paper.

From a broader perspective the Kummer construction furnishes the prototypical example of the appearance of orbifold singularities in non-collapsing sequences of Einstein $4$--manifolds. In \cite[Theorem C]{Anderson:Ricci:Bounds} Anderson showed that a sequence of Einstein $4$--manifolds $(M_i,g_i)$ with a uniform lower bound on volume and upper bounds on diameter and Euler characteristic converges (up to subsequences) to an Einstein $4$--orbifold $M_\infty$ with finitely many singular points. The formation of orbifold singularities is modelled on complete Ricci-flat ALE spaces which appear as rescaled limits, or ``bubbles'', of the sequence $(M_i,g_i)$ around points that approach one of the singularities of the orbifold $M_\infty$.

In the Ricci-flat case collapsing can also occur. Anderson \cite[Theorem II]{Anderson:L2:structure} showed that every sequence of Ricci-flat metrics $(M,g_i)$ of unit volume but unbounded diameter collapses everywhere, \ie $\text{inj}_{g_i}(x)\ra 0$ for all $x\in M$. The collapse is in the sense of Cheeger--Gromov outside finitely many points $x_1,\dots, x_n$, \ie $\text{inj}_{g_i}(x)\ra 0$ and $\text{inj}_{g_i}(x)^2|\text{Rm}_{g_i}|_{g_i}(x)\leq \epsilon_0$ for all $x\in M\setminus \{ x_1,\dots, x_n\}$, for a universal constant $\epsilon_0>0$. In fact, Cheeger and Tian \cite[Theorems 0.1 and 0.8]{Cheeger:Tian} have shown that the collapse occurs with \emph{bounded} curvature away from a definite number of points.

Contrary to the case of orbifold singularities, almost nothing is known about the structure of the singular points arising in collapsing sequences of Ricci-flat metrics \cite[\S 6]{Anderson:Survey}. One would expect that the geometry around these points is modelled on complete Ricci-flat manifolds with non-maximal volume growth. A simple example of the expected  phenomena was suggested by Page \cite{Page:D2:ALF} in 1981. Consider the Kummer construction of Ricci-flat metrics on the K3 surface along a family of split tori $T^4 = T^3 \times S^1_\ell$ with a circle factor of length $\ell\ra 0$. We can then think of the $2$--spheres arising in the resolution of the $16$ singularities of $T^4/\Z_2$ as coming in pairs aligned along the collapsing circle over each of the $8$ singular points of $T^3/\Z_2$. If we now rescale the sequence of K\"ahler Ricci-flat metrics on the K3 surface by $\ell^{-2}$ around one of these pairs, Page suggests, in the limit $\ell\ra 0$ we should obtain a complete Ricci-flat (hyperk\"ahler) metric on a noncompact space which at infinity looks like $(\R^3\times S^1)/\Z_2$.

In this paper we regard this example as a simple case of a more general construction of sequences of Ricci-flat metrics on the K3 surface that collapse to a $3$--dimensional limit. In this more general construction, the $4$--torus $T^3\times S^1_\ell$ is replaced by a non-trivial circle bundle over a (punctured) $3$--torus and the role of Page's ``periodic but nonstationary gravitational instanton'' is played by other ALF gravitational instantons.

A gravitational instanton is a complete hyperk\"ahler $4$--manifold with decaying curvature at infinity. Since every hyperk\"ahler manifold is in particular Ricci-flat, gravitational instantons have constrained volume growth: the volume of a geodesic ball of radius $r$ grows at most as $r^4$. Gravitational instantons of maximal volume growth are the ALE spaces constructed and classified by Kronheimer \cite{Kronheimer} following earlier work of Eguchi--Hanson, Gibbons--Hawking and Hitchin. We have seen how ALE spaces arise as models for the formation of orbifold singularities of non-collapsed sequences of Einstein $4$--manifolds. In \cite{Minerbe:ALF} Minerbe showed that if $\vol \big( B_r(p)\big)=O(r^a)$ for some $3\leq a <4$ and all $p$, then we must have $a=3$. Gravitational instantons of cubic volume growth are called ALF. By \cite{Minerbe:ALF} the (unique) end of an ALF space looks like a circle fibration over the complement of a ball in $\R^3$ or $\R^3/\Z_2$ with fibres of asymptotically finite length. If the base of the circle fibration at infinity is $\R^3$ (respectively, $\R^3/\Z_2$) then we say that the ALF space is of cyclic (dihedral) type, since the boundary of large geodesic balls is diffeomorphic to $S^3/\Gamma$, where $\Gamma \subset SU(2)$ is a cyclic group in the first case and a binary dihedral group in the second.

The prototypical example of an ALF space of cyclic type is the Taub--NUT metric on $\R^4$. This metric is explicit and the circle fibration at infinity is induced by the Hopf projection $S^3 \ra S^2$. The first example of an ALF metric of dihedral type was found by Atiyah--Hitchin \cite{Atiyah:Hitchin} by studying moduli spaces of magnetic monopoles on $\R^3$, \ie the solutions of the dimensional reduction of the Yang--Mills self-duality equations from $4$ to $3$ dimensions. The Atiyah--Hitchin manifold is diffeomorphic to the complement of a Veronese $\R\PP^2$ in $\Sph^4$ and the metric is explicitly given in terms of elliptic integrals.  

Recently ALF gravitational instantons have been the focus of intense research with the aim of constructing and classifying examples. Minerbe \cite{Minerbe:Ak} classified ALF spaces of cyclic type. These are all explicitly given by the Gibbons--Hawking construction of hyperk\"ahler $4$--manifolds with a triholomorphic circle action \cite{Gibbons:Hawking}. Most (if not all) the known methods of construction of hyperk\"ahler metrics have been applied to the dihedral ALF case: twistor methods \cite{Hitchin:Einstein,Cherkis:Hitchin}, hyperk\"ahler quotient constructions \cite{Dancer}, gauge-theoretic constructions as in the case of the Atiyah--Hitchin manifold \cite{Cherkis:Hitchin}, Kummer-type constructions \cite{Biquard:Minerbe} and complex Monge--Amp\`ere methods \cite{Auvray:I,Auvray:II}. For example, Page's ``periodic but nonstationary'' gravitational instantons of \cite{Page:D2:ALF}, more commonly known as $D_2$ ALF spaces, were first constructed rigorously by Hitchin \cite{Hitchin:Einstein} using twistor methods and more recently by Biquard--Minerbe \cite{Biquard:Minerbe} using an extension of the Kummer construction to non-compact spaces. The $D_2$ ALF spaces can also be thought of as the moduli spaces of centred charge $2$ $SO(3)$ monopoles on $\R^3$ with two singularities endowed with their natural $L^2$--metric \cite{Cherkis:Kapustin}. The fact that all these constructions yield equivalent families of ALF metrics was shown only recently by Chen--Chen \cite{Chen:Chen:II}. 

Despite this rich theory of ALF gravitational instantons, until now it has remained unclear how they can appear as models for the formation of singularities in collapsing sequences of hyperk\"ahler metrics on the K3 surface. The aim of this paper is to exploit singular perturbation methods to construct examples of Ricci-flat metrics on the K3 surface collapsing to a $3$--dimensional limit and exhibit ALF gravitational instantons as the ``bubbles'' appearing in the process.

\begin{theorem}\label{thm:Main}
Every collection of $8$ ALF spaces of dihedral type $M_1,\dots, M_8$ and $n \leq 16$ ALF spaces of cyclic type $N_1,\dots, N_n$ satisfying
\[
\sum_{j=1}^8{\chi(M_j)} + \sum_{i=1}^n{\chi (N_i)} =24
\]
arises as the collection of ``bubbles'' forming in a sequence of hyperk\"ahler metrics on the K3 surface which collapse to $T^3/\Z_2$ with bounded curvature away from $n+8$ points.
\end{theorem}
We refer to Theorem \ref{thm:Main:precise} for a more precise statement.

In \cite{Gross:Wilson} Gross--Wilson studied hyperk\"ahler metrics on elliptic K3 surfaces with fibres of small size. They considered the generic case when all singular fibres ($24$ of them) are of Kodaira type $I_1$ (\ie a pinched torus). The hyperk\"ahler metric is approximated by a semi-flat metric on the locus of the smooth fibres and by a certain (incomplete) explicit hyperk\"ahler metric, the Ooguri--Vafa metric, in the neighbourhood of each singular fibre. As the size of the fibres converges to zero, the K3 surface collapses to a metric on $S^2$ (the base of the elliptic fibration) with $24$ singular points. To the knowledge of the author, besides Gross--Wilson's work, Theorem \ref{thm:Main} is the only study of collapsing sequences of hyperk\"ahler metrics on the K3 surface.

Now, one way to make precise Page's observations in \cite{Page:D2:ALF} about the Kummer construction for a degenerating family of tori is to consider a gluing construction in which one glues $8$ copies of the $D_2$ ALF space to the $\Z_2$ quotient of the trivial circle bundle $T^3 \times S^1$ over the flat $3$--torus. The proof of Theorem \ref{thm:Main} is also based on a gluing construction. In order to allow for more general ALF spaces to appear as rescaled limits, the main idea is to replace $T^3\times S^1$ with an (incomplete) background hyperk\"ahler metric on a \emph{non-trivial} circle bundle over a punctured $3$--torus. The tool to construct such a background metric is the Gibbons--Hawking construction of hyperk\"ahler metrics with a triholomorphic $S^1$ symmetry, \ie an isometric circle action that also preserves the $2$--sphere of complex structures compatible with the metric. The $S^1$--invariant hyperk\"ahler metrics we seek are explicitly given in terms of a positive harmonic function $h$ on $T^3$ with prescribed singularities at a finite number of points. For most configurations of punctures the harmonic function $h$ becomes negative somewhere. However, by multiplying $h$ by a small number $\epsilon>0$ (which geometrically corresponds to making the circle fibres have small length) it is possible to construct highly collapsed hyperk\"ahler metrics $g^\textup{gh}_\epsilon$ outside of an arbitrarily small neighbourhood of the punctures. Furthermore, the construction of this background metric can be made invariant under the action of an involution. 

The key observation now is that the asymptotic model of any ALF metric (up to a double cover in the dihedral case) can be written in Gibbons--Hawking coordinates. By choosing the configuration of punctures appropriately it is then possible to glue in copies of ALF spaces to extend the Gibbons--Hawking metric $g^\textup{gh}_\epsilon$ to an approximately hyperk\"ahler metric $g_\epsilon$: close to a fixed point of the $\Z_2$--action on $T^3$ we glue in an ALF space of dihedral type (this explains why we need $8$ of them in Theorem \ref{thm:Main}); close to a puncture which is not a fixed point of the $\Z_2$--action we glue in an ALF space of cyclic type. The Euler characteristic constraint in the statement of Theorem \ref{thm:Main} is necessary for the resulting $4$--manifold to have the same Euler characteristic as the K3 surface, but it can also be reinterpreted as the necessary and sufficient condition for the existence of the harmonic function $h$ in the first place. 

The approximate solution $g_\epsilon$ is then deformed into an exact hyperk\"ahler metric by means of the Implicit Function Theorem. Since some of the ALF spaces are not biholomorphic to their asymptotic model outside a compact set, it is necessary to set up the problem as a gluing problem for hyperk\"ahler structures, rather than the most standard procedure (as in the classical Kummer construction) of first constructing a complex surface using complex geometry and then solving a complex Monge--Amp\`ere equation on this given complex manifold.

\begin{remark*}
At least in some form this ``Gibbons--Hawking approximation'' of hyperk\"ahler metrics on the K3 surface seems to be known to physicists in the context of the duality between M theory compactified on the K3 surface and Type IIA String theory compactified on $T^3/\Z_2$. For example, in \cite{Sen} Sen discusses the physical interpretation of dihedral ALF spaces thought of as a ``superposition'' of Taub--NUT spaces and the Atiyah--Hitchin manifold, \cf Remark \ref{rmk:Dk}.
\end{remark*}

\begin{remark*}
One can also wonder what happens when we start from an arbitrary orientable flat $3$--manifold instead of a $3$--torus. There are $6$ of these: in the notation of \cite[\S 3.5]{Wolf} they are $\mathcal{G}_1=T^3$, $\mathcal{G}_i=T^3/\Z_i$ for $i=2,3,4$, $\mathcal{G}_5=T^3/\Z_6$ and $\mathcal{G}_6=T^3/(\Z_2\times \Z_2)$. Only $\mathcal{G}_1$ has $b_1 =3$, $b_1 (\mathcal{G}_i)=1$ in all other cases except for $\mathcal{G}_6$ which has purely torsion first homology \cite[Equation (2.5)]{Luft:Sjerve}. By working on the $3$--torus equivariantly with respect to a finite group action, the Gibbons--Hawking construction then yields (incomplete) Ricci-flat metrics on circle bundles over a punctured flat $3$--manifold $M$ which are hyperk\"ahler only when $M=T^3$, K\"ahler if $M=\mathcal{G}_i$ for $i=2,3,4,5$ and have generic holonomy when $M=\mathcal{G}_6$. Moreover, Luft--Sjerve \cite[Theorem 1.1]{Luft:Sjerve} have shown that only $\mathcal{G}_1, \mathcal{G}_2$ and $\mathcal{G}_6$ admit an involution with finitely many fixed points ($8$, $4$ and $2$ of them, respectively). Hence only in these $3$ cases are we able to construct background Ricci-flat metrics that can be extended to complete metrics by gluing in copies of ALF spaces of cyclic and dihedral type. On the other hand, Hitchin \cite[Theorem 1]{Hitchin:Compact:Einstein} showed that the only Ricci-flat $4$--manifolds covered by the K3 surface are the Enriques surfaces (quotients of a K3 surface by an involution without fixed points) with their K\"ahler Ricci-flat metrics and the quotient of an Enriques surface by an anti-holomorphic involution without fixed points. Carrying out our gluing construction equivariantly with respect to a finite group action then allows us to produce collapsing sequences of Ricci-flat metrics on an Enriques surface (the metrics are K\"ahler in this case) and its quotient by an anti-holomorphic involution: the collapsed limit is $\mathcal{G}_2/\Z_2$ and $\mathcal{G}_6/\Z_2$, respectively, and ALF gravitational instantons appear as ``bubbles''.
\end{remark*}

We leave aside for future work the question of understanding the relation between the metric degenerations described in this paper and degenerations of a compatible complex structure on the K3 surface. Similarly, it would be very interesting to understand to what extent the collapsing behaviour exhibited in this paper is typical of an arbitrary sequence of Ricci-flat metrics on the K3 surface collapsing to a $3$--dimensional limit.

We give instead an application of our gluing construction to the theory of minimal surfaces and harmonic maps. It is well known that holomorphic submanifolds of a K\"ahler manifold minimise volume in their homology class. A classical problem in minimal surface theory is to understand to what extent area minimising surfaces (and more generally stable minimal surfaces) in K\"ahler manifolds must be (anti)holomorphic. For example, in 1993 Yau asked whether it is possible to classify all stable minimal $2$--spheres in a simply connected K\"ahler Ricci-flat manifold \cite[Question 64]{Yau}. In \cite{Micallef} Micallef showed that every stable minimal surface in a flat $4$--torus must be holomorphic for some complex structure compatible with the metric. For some time there was hope to prove a similar result in the case of the K3 surface endowed with a hyperk\"ahler metric. Eventually, Micallef--Wolfson \cite{Micallef:Wolfson:2} showed that this is not the case. A simple application of our gluing construction allows us to give an alternative (simpler) counterexample: there exist hyperk\"ahler metrics on the K3 surface that admit a strictly stable minimal sphere which is not holomorphic with respect to any complex structure compatible with the metric, \cf Theorem \ref{thm:Minimal:Surfaces}.

In fact each of these minimal spheres is the image of a smooth harmonic map $u\co S^2 \ra \textup{K3}$ satisfying
\begin{equation}\label{eqn:Quaternionic:minimal:surfaces}
du\circ J_{S^2} = -\left( x_1\, J_1\circ du + x_2\, J_2\circ du + x_3\, J_3\circ du \right),
\end{equation}
where $(x_1,x_2,x_3)$ defines the standard embedding of $S^2$ as the unit sphere in $\R^3$ and $(J_1,J_2,J_3)$ is the triple of parallel complex structures associated to a hyperk\"ahler structure on the K3 surface. A map $u$ satisfying \eqref{eqn:Quaternionic:minimal:surfaces} defines a homogeneous triholomorphic map $F(x,x_4)=u(x/|x|)$ from the quaternions $\HH$ into the K3 surface with singular sets $\{ x=0\}$ of codimension $3$. Here a triholomorphic map is a harmonic map between hyperk\"ahler manifolds satisfying a quaternionic del-bar equation, \cf \S \ref{sec:Minimal:surfaces} for a precise definition. Triholomorphic maps (and their dimensional reduction to Fueter maps from $3$--manifolds into hyperk\"ahler manifolds) arise in different contexts related to gauge theory in higher dimensions, \cf for example \cite{OFarrill} and \cite[\S 6]{Donaldson:Segal}, the definition of enumerative invariants of hyperk\"ahler manifolds \cite{Tian:GW:inv:HK} and the hyperk\"ahler Floer theory of \cite{Salamon:al}. For these geometric applications it is essential to understand the regularity theory of limits of smooth triholomorphic maps. The importance of the existence or non-existence of homogeneous triholomorphic maps with a singular set of codimension $3$ was stressed in the recent work of Bellettini--Tian \cite[Remark 1.2]{Bellettini:Tian}. In Theorem \ref{thm:Triholomorphic} we construct the first known examples of such singular triholomorphic maps with compact target.

\subsection*{Plan of the paper}

As we have already mentioned, in this paper we will need to glue hyperk\"ahler structures rather than solving a complex Monge--Amp\`ere equation on a given complex manifold. In Section \ref{sec:Definite:triples}, following Donaldson \cite{Donaldson:2:forms}, we explain how to set up the problem of deforming approximately hyperk\"ahler metrics based on the notion of definite triples.

Section \ref{sec:ALF} is a detailed summary of the theory of ALF spaces: we give precise definitions, describe detailed asymptotics for such metrics and recall the construction and classification of examples.

In Section \ref{sec:GH:3:torus} we use the Gibbons--Hawking ansatz to construct (incomplete) hyperk\"ahler metrics on circle bundles over a punctured $3$--torus. In Section \ref{sec:Approximate:solution} we use ALF spaces of cyclic and dihedral type together with the metrics constructed in Section \ref{sec:GH:3:torus} to produce families of approximately hyperk\"ahler metrics. In Section \ref{sec:Deformation} we use analysis to deform these approximate solutions into exact hyperk\"ahler metrics. This is done by means of an Implicit Function Theorem in weighted H\"older spaces. As usual in gluing problems, most of the work goes into showing that the relevant linear operator has no small eigenvalues as $\epsilon\ra 0$ and the geometry degenerates. 

Finally, Section \ref{sec:Minimal:surfaces} contains the proof of Theorem \ref{thm:Minimal:Surfaces} about the existence of non-holomorphic strictly stable minimal spheres and of Theorem \ref{thm:Triholomorphic} about the existence of homogeneous triholomorphic maps with a singular set of codimension $3$.

\subsection*{Acknowledgements} The author wishes to thank Bobby Acharya, Mark Haskins and Johannes Nordstr\"om for an inspiring conversation at the Mathematisches Forschungsinstitut Oberwolfach in February 2015 which was the original inspiration for this work. He also wishes to thank Mark Haskins for reading an earlier version of the paper and for many discussions and suggestions for improvement. Discussions with Mark Haskins on this work and related topics were also made possible thanks to the support of his EPSRC grant EP/L001527/1, ``Singular spaces of special and exceptional holonomy''. Finally, the author thanks Thomas Walpuski for mentioning to him the relevance of the minimal sphere in the double-cover of the Atiyah--Hitchin manifold for the regularity theory of Fueter and triholomorphic maps. The paper is based on work supported by the National Science Foundation under Grant No. DMS-1440140 while the author was in residence at the Mathematical Sciences Research Institute in Berkeley, California, during the Spring 2016 semester, and Grant No. DMS-1608143.

\section{Definite triples and hyperk\"ahler structures on $4$--manifolds}\label{sec:Definite:triples}

The standard approach in the Kummer construction of K\"ahler Ricci-flat metrics on the K3 surface is to proceed in two steps. First one constructs a complex surface with vanishing first Chern class by blowing up the singularities of a flat orbifold $T^4/\Z_2$. On this given complex manifold one then constructs a K\"ahler Ricci-flat metric by solving a complex Monge-Amp\`ere equation. Some of the building blocks we are going to use in the gluing construction of this paper are not biholomorphic to their asymptotic model outside a compact set. This will force us to adopt a different strategy and glue hyperk\"ahler structures all together. In this initial preliminary section we explain how Donaldson \cite{Donaldson:2:forms} suggested an approach to this problem, based on the notion of definite triples.

Recall that the space of $2$--forms on an oriented $4$--dimensional vector space carries a natural non-degenerate bilinear form of signature $(3,3)$.

\begin{definition}\label{def:Definite:Triple}
Let $(M^4,\mu_0)$ be an oriented $4$--manifold with volume form $\mu_0$. A \emph{definite triple} is a triple $\triple{\omega}=(\omega_1,\omega_2,\omega_3)$ of $2$--forms on $M$ such that $\text{Span}(\omega_1,\omega_2,\omega_3)$ is a $3$--dimensional positive definite subspace of $\Lambda^2 T^\ast_x M$ at every point $x\in M$.
\end{definition}

Given a triple $\triple{\omega}$ of $2$--forms on $(M,\mu_0)$ we consider the matrix $Q\in \Gamma \big( M,\text{Sym}^2(\R^3) \big)$ defined by
\begin{equation}\label{eq:Intersection:Matrix}
\tfrac{1}{2}\,\omega_i \wedge \omega_j = Q_{ij}\,\mu_0.
\end{equation}
$\triple{\omega}$ is a definite triple if and only if $Q$ is a positive definite matrix.

To every definite triple $\triple{\omega}$ we associate a volume form $\mu_{\triple{\omega}}$ by
\begin{equation}\label{eq:Volume:form}
\mu_{\triple{\omega}} = \left( \det Q \right)^{\frac{1}{3}}\mu_0
\end{equation}
and the new matrix $Q_{\triple{\omega}}=\left( \det{Q} \right) ^{-\frac{1}{3}}Q$ which satisfies \eqref{eq:Intersection:Matrix} with $\mu_{\triple{\omega}}$ in place of $\mu_0$. Note that the volume form $\mu_{\triple{\omega}}$ and the matrix $Q_{\triple{\omega}}$ are independent of the choice of volume form $\mu_0$. We refer to $\mu_{\triple{\omega}}$ and $Q_{\triple{\omega}}$ as the \emph{associated volume form} and \emph{intersection matrix} of the definite triple $\triple{\omega}$.

Now, let $(M^4,\mu_0)$ be an oriented $4$--dimensional manifold. It is well known that the choice of a $3$--dimensional positive definite subspace of $\Lambda^2 T^\ast _x M$ for all $x\in M$ is equivalent to the choice of a conformal class on $M$. Thus every definite triple defines a Riemannian metric $g_{\triple{\omega}}$ by requiring that $\text{Span}(\omega_1,\omega_2,\omega_3)|_x = \Lambda^+T^\ast_x M$ for all $x\in M$ and $\dvol_{g_{\triple{\omega}}}=\mu_{\triple{\omega}}$.  

\begin{definition}\label{eq:Hyperkahler}
A definite triple $\triple{\omega}$ is said
\begin{enumerate}
\item \emph{closed} if $d\omega_i=0$ for $i=1,2,3$;
\item an \emph{$SU(2)$--structure} if $Q_{\triple{\omega}}\equiv \text{id}$;
\item \emph{hyperk\"ahler} if it is both closed and an $SU(2)$--structure.
\end{enumerate}
\end{definition}

The metric $g_{\triple{\omega}}$ associated to a hyperk\"ahler triple is hyperk\"ahler, in the sense that it has holonomy contained in $Sp(1)\simeq SU(2)$.

\subsection{The deformation problem}

In Section \ref{sec:Approximate:solution} we will construct closed definite triples $\triple{\omega}$ which are approximately hyperk\"ahler, in the sense that the intersection matrix $Q_{\triple{\omega}}$ is close to the identity. We now explain how to formulate the problem of deforming such a triple $\triple{\omega}$ into a hyperk\"ahler structure. 

Let $\triple{\omega}$ be a closed definite triple on a $4$--manifold $M$ and assume that $\| Q_{\triple{\omega}}-\text{id}\| _{C^0} < \sigma$ for some small $\sigma>0$. We want to deform $\triple{\omega}$ into a hyperk\"ahler triple, \ie we look for a triple of closed $2$--forms $\triple{\eta}=(\eta_1,\eta_2,\eta_3)$ on $M$ such that
\begin{equation}\label{eq:Nonlinear:1}
\tfrac{1}{2}\left( \omega_i + \eta_i\right) \wedge \left( \omega_j + \eta_j\right) = \delta_{ij}\, \mu_{\triple{\omega}}.
\end{equation}

Decompose $\triple{\eta}$ into self-dual and anti-self dual parts $\triple{\eta}=\triple{\eta}^+ + \triple{\eta}^-$ with respect to $g_{\triple{\omega}}$. The self-dual part can be written in terms of a $M_{3\times 3}(\R)$--valued function $A$ by
\[
\eta_i^+=\sum_{j=1}^3{A_{ij}\, \omega_j}.
\]

Denote by $\triple{\eta}^-\ast\triple{\eta}^-$ the symmetric $(3\times 3)$--matrix with entries $(\tfrac{1}{2}\,\eta_i^-\wedge\eta_j^-)/\mu_{\triple{\omega}}$. Then we can rewrite \eqref{eq:Nonlinear:1} as
\begin{equation}\label{eq:Nonlinear:2}
Q_{\triple{\omega}}+Q_{\triple{\omega}}\, A^T+A\, Q_{\triple{\omega}}+A\, Q_{\triple{\omega}}\,A^T+\triple{\eta}^-\ast\triple{\eta}^-=\text{id}.
\end{equation}

Now, consider the map
\[
M_{3\times 3}(\R)\longrightarrow Sym^{2}(\R^3);\qquad  A\longmapsto Q_{\triple{\omega}}\,A^T+A\, Q_{\triple{\omega}}+A\, Q_{\triple{\omega}}\, A^T
\]
and its differential $A\mapsto Q_{\triple{\omega}}\, A^T+A\, Q_{\triple{\omega}}$. Since $Q_{\triple{\omega}}$ is arbitrarily close to the identity, this linear map induces an isomorphism $Sym^{2}(\R^3) \ra Sym^{2}(\R^3)$ for $\sigma$ sufficiently small. We can therefore define a smooth function $\mathcal{F}\co Sym^{2}(\R^3) \ra Sym^{2}(\R^3)$ such that $Q_{\triple{\omega}}\,A^T+A\, Q_{\triple{\omega}} +A\, Q_{\triple{\omega}}\, A^T=S$ if and only if $A=\mathcal{F}(S)$.

\begin{remark*}
When $\triple{\omega}$ is hyperk\"ahler (thus $Q_{\triple{\omega}}=\text{id}$) the kernel of $A\mapsto Q_{\triple{\omega}}\, A^T+A\, Q_{\triple{\omega}}$ corresponds to infinitesimal hyperk\"ahler rotations.
\end{remark*}

Hence we reformulate \eqref{eq:Nonlinear:2} as
\begin{equation}\label{eq:Nonlinear:3}
\triple{\eta}^+=\mathcal{F}\left( (\text{id}-Q_{\triple{\omega}})-\triple{\eta}^-\ast\triple{\eta}^-\right).
\end{equation}

Now, let $\mathcal{H}^+_{\triple{\omega}}$ be the space of self-dual harmonic $2$--forms with respect to $g_{\triple{\omega}}$. If a solution of \eqref{eq:Nonlinear:2} exists on a compact manifold $M$ then $M$ must be either a $4$--torus or a K3 surface with the standard orientation and therefore $\mathcal{H}^+_{\triple{\omega}}$ is $3$--dimensional. Since $\omega_1,\omega_2,\omega_3$ are closed and self-dual (therefore harmonic) and linearly independent (since $\triple{\omega}$ is a definite triple) we deduce that $\mathcal{H}^+_{\triple{\omega}}$ consist of constant linear combinations of $\omega_1,\omega_2,\omega_3$.

By Hodge theory with respect to $g_{\triple{\omega}}$ we can finally rewrite \eqref{eq:Nonlinear:3} as the \emph{elliptic} equation
\begin{equation}\label{eq:Nonlinear}
d^+\triple{a}+\triple{\zeta}=\mathcal{F}\left( (\text{id}-Q_{\triple{\omega}})-\triple{\eta}^-\ast\triple{\eta}^-\right), \qquad d^\ast\triple{a}=0,
\end{equation}
for a triple $\triple{a}$ of $1$--forms on $M$ and a triple $\triple{\zeta} \in \mathcal{H}^+_{\triple{\omega}}\otimes \R^3$. Here $2\, d^+a = da+\ast da$ is the self-dual part of $da$.

\begin{remark*}
Note that in general it is necessary to deform the cohomology classes of $\omega_1,\omega_2,\omega_3$ since every hyperk\"ahler triple must satisfy
\[
\tfrac{1}{2}\langle\, [\omega_i]\cup [\omega_j],[M]\,\rangle =\delta_{ij} \vol_{g_{\triple{\omega}}} (M).
\]
\end{remark*}

The linearisation of \eqref{eq:Nonlinear} is
\begin{equation}\label{eq:Linearisation}
(D\oplus \text{id})\otimes \R^3 \co \left( \Omega^1(M)\oplus \mathcal{H}^+_{\triple{\omega}}\right) \otimes \R^3 \longrightarrow \left( \Omega^0(M)\oplus \Omega^+(M) \right) \otimes \R^3,
\end{equation}
where $D$ is the Dirac-type operator
\begin{equation}\label{eq:Dirac}
D=d^\ast \oplus d^+\co \Omega^1(M)\longrightarrow \Omega^0(M)\oplus \Omega^+(M).
\end{equation}
Note that the operator in \eqref{eq:Linearisation} is always surjective with kernel consisting of harmonic $1$--forms.

\section{ALF gravitational instantons}\label{sec:ALF}

In this section we collect known results about gravitational instantons of type ALF, with an emphasis on their asymptotic geometry. ALF gravitational instantons will appear as local models for the geometry of high curvature regions in sequences of hyperk\"ahler metrics on K3 collapsing to a $3$--dimensional limit.

A \emph{gravitational instanton} is a complete hyperk\"ahler $4$--manifold $(M,g)$ with decaying Riemannian curvature at infinity. The minimum requirement (automatically satisfied for rescaled limits of Einstein metrics on $4$--manifolds with bounded Euler characteristic by the Chern--Gauss--Bonnet formula) is that $(M,g)$ has finite \emph{energy} $\| \text{Rm}\|_{L^2}$. In order to say something about the structure of gravitational instantons it has often been necessary to strengthen this finite energy assumption to \emph{faster than quadratic curvature decay} $|\text{Rm}|=O(r^{-2-\epsilon})$, $\epsilon>0$ (or a slightly weaker finite weighted energy assumption). Note however that there are examples of gravitational instantons which do \emph{not} satisfy this stronger decay assumption \cite[Theorem 1.5]{Hein}.

Since hyperk\"ahler manifolds are in particular Ricci-flat, gravitational instantons have only one end and constrained \emph{volume growth}: the volume of a geodesic ball of radius $r$ can grow at most as $r^4$ and at least linearly. An initial rough classification of gravitational instantons can be given in terms of their volume growth. The gravitational instantons of maximal volume growth are the ALE spaces classified by Kronheimer \cite{Kronheimer} following earlier work of Eguchi--Hanson, Gibbons--Hawking and Hitchin. Under the assumption of faster than quadratic curvature decay (or a slightly weaker finite weighted energy assumption) Minerbe \cite[Theorem 0.1]{Minerbe:ALF} has shown that if we assume $\vol \big( B_r(p)\big)=O(r^a)$ for some $3\leq a<4$ and all $p$, then $a=3$. Minerbe also described the asymptotic geometry of gravitational instantons of cubic volume growth and faster than quadratic curvature decay: they are all ALF spaces, in the following sense.

\begin{definition}\label{def:ALF}
A gravitational instanton $(M,g)$ is called \emph{ALF} if there exists a compact set $K\subset M$, $R>0$ and a finite group $\Gamma < O(3)$ acting freely on $\Sph^2$ such that $M\setminus K$ is the total space of a circle fibration $\pi\co M\setminus K \ra (\R^3\setminus B_R)/\Gamma$ and the metric is asymptotically a Riemannian submersion
\begin{equation}\label{eq:ALF:weak}
g=\pi^\ast g_{\R^3/\Gamma} + \theta^2 + O(r^{-\tau})
\end{equation}
for a connection $\theta$ on $\pi$ and some $\tau >0$. There are two possibilities for the finite group $\Gamma$: if $\Gamma=\text{id}$ we say that $M$ is an ALF gravitational instanton of \emph{cyclic type}; if $\Gamma = \Z_2$ we say that $M$ is an ALF gravitational instanton of \emph{dihedral type}.
\end{definition}

We are interested in refining these asymptotics. In order to describe a more precise model for the end of an ALF gravitational instanton it is necessary to recall the explicit construction of $4$--dimensional hyperk\"ahler metrics with a triholomorphic circle action known as the Gibbons--Hawking ansatz. We will use this same ansatz later in the paper to construct (incomplete) hyperk\"ahler metrics on circle bundles over a punctured $3$--torus.

\subsection{The Gibbons--Hawking ansatz}

The Gibbons--Hawking ansatz describes $4$--dimensional hyperk\"ahler metrics with an isometric $S^1$--action that also preserves the whole hyperk\"ahler structure. Such an $S^1$ action is therefore called \emph{triholomorphic}.

Let $U$ be an open set of $\R^3$ and $\pi\co P\ra U$ be a principal $U(1)$--bundle. Suppose that there exists a positive harmonic function $h$ on $U$ such that $\ast dh$ is the curvature $d\theta$ of a connection $\theta$ on $P$. Then
\begin{subequations}\label{eq:Gibbons:Hawking}
\begin{equation}\label{eq:Gibbons:Hawking:metric}
g^{\textup{gh}}= h\, \pi^\ast g_{\R^3} + h^{-1}\theta^2
\end{equation}
is a hyperk\"ahler metric. Indeed, we can exhibit an explicit hyperk\"ahler triple $\triple{\omega}^{\textup{gh}}$ that induces the metric $g^{\textup{gh}}$. Fix coordinates $(x_1,x_2,x_3)$ on $U\subset \R^3$ and define
\begin{equation}\label{eq:Gibbons:Hawking:triple}
\omega^{\textup{gh}}_i = dx_i \wedge \theta + h\, dx_j \wedge dx_k.
\end{equation}
\end{subequations}
Here and in the rest of the paper we use the convention that for every $i=1,2,3$ the indices $j,k$ are chosen so that $\epsilon_{ijk}=1$. One can check explicitly that $\triple{\omega}^{gh}$ defines an $SU(2)$--structure and it induces the Riemannian metric $g^{gh}$. Moreover, the requirement that $\triple{\omega}^{gh}$ is also closed is equivalent to the abelian \emph{monopole equation}
\begin{equation}\label{eq:Monopole}
\ast dh = d\theta
\end{equation}
The fibre-wise circle action on $P$ preserves $\triple{\omega}^\textup{gh}$ and $\pi$ is nothing but a hyperk\"ahler moment map for this action. Conversely, every $4$--dimensional hyperk\"ahler metric with a triholomorphic circle action is described by \eqref{eq:Gibbons:Hawking}.

The basic example of the Gibbons--Hawking construction is given in terms of so-called Dirac monopoles on $\R^3$. Fix a set of distinct points $p_1,\dots,p_n$ in $\R^3$ and consider the harmonic function
\[
h=\lambda + \sum_{j=1}^n{\frac{k_j}{2|x-p_j|}},
\]
where $\lambda>0$ and $k_1,\dots, k_n$ are constants. Since $\R^3\setminus \{ p_1,\dots, p_n\}$ has non-trivial second homology, we must require $k_j\in\Z$ for all $j$ in order to be able to solve \eqref{eq:Monopole}. If these integrality constraints are satisfied then $\ast dh$ defines the curvature $d\theta$ of a connection $\theta$ (unique up to gauge transformations) on a principal $U(1)$--bundle $P$ over $\R^3\setminus \{ p_1,\dots, p_n\}$ which restricts to the principal $U(1)$--bundle associated with the line bundle $\mathcal{O}(k_j)\ra S^2$ on a small punctured neighbourhood of $p_j$. The pair $(h,\theta)$ is a solution of \eqref{eq:Monopole} which we call a \emph{Dirac monopole} with singularities at $p_1,\dots, p_n$.

The Gibbons--Hawking ansatz \eqref{eq:Gibbons:Hawking} associates a hyperk\"ahler metric $g^{\textup{gh}}$ to every Dirac monopole on the open set where $h>0$. When $k_j > 0$ then $g^{\textup{gh}}$ is certainly defined on the restriction of $P$ to a small punctured neighbourhood of $p_j$. By a change of variables one can check that $g^{gh}$ can be extended to a smooth (orbifold) metric modelled on $\C^2/\Z_{k_j}$ by adding a single point. In particular $g^{\textup{gh}}$ is a complete metric whenever $\lambda\geq 0$ and $k_j=1$ for all $j=1,\dots,n$. One can check that $g^{\textup{gh}}$ is an ALE metric when $\lambda=0$ and an ALF metric of cyclic type when $\lambda>0$. Note also that when $\lambda >0$ we can always rescale the metric so that $\lambda=1$.

\subsection{Families of ALF gravitational instantons}

We are now going to use the Gibbons--Hawking ansatz to define a refined asymptotic model for ALF gravitational instantons.

Let $H^k$ be the total space of the principal $U(1)$--bundle associated with the line bundle $\mathcal{O}(k)$ over $S^2$ radially extended to $\R^3\setminus B_R$ for any $R>0$. $\theta_k$ will denote the $SO(3)$--invariant connection on $H^k$. The Gibbons--Hawking ansatz \eqref{eq:Gibbons:Hawking} yields a hyperk\"ahler metric
\begin{equation}\label{eq:Metric:Infinity}
g_k = \left( 1+\frac{k}{2\rho} \right) (d\rho^2 + \rho^2 g_{\Sph^2}) + \left( 1+\frac{k}{2\rho} \right)^{-1} \theta_k^2
\end{equation}
on $H^k$ for all $k\in \Z$. Here $\rho$ is a radial function on $\R^3$. We denote by $\triple{\omega}_k$ the associated hyperk\"ahler triple defined by \eqref{eq:Gibbons:Hawking:triple}. Note that we could replace the harmonic function $1+\frac{k}{2\rho}$ with $\lambda+\frac{k}{2\rho}$ for any $\lambda>0$ but we can always reduce to the case $\lambda=1$ by scaling.

Finally, on $H^{2k}$ we consider the $\Z_2$--action which is defined as the simultaneous standard involutions on the base $\R^3$ and the fibre. Here the involution on the fibre $S^1 =\R /2\pi\Z$ is the one induced by the standard involution on the universal cover $\R$. Throughout the paper we refer to this as the standard involution of $S^1$. 

\begin{definition}\label{def:Ak:Dk}
Let $(M^4,g)$ be an ALF gravitational instanton of cyclic type. By scaling assume that the length of the circle fibres at infinity is $1$.
\begin{enumerate}
\item We say that $M$ is of type $A_k$ for some $k\geq -1$ if there exists a compact set $K\subset M$, $R>0$ and a diffeomorphism $\phi\co H^{k+1} \ra M\setminus K$ such that
\[
|\nabla ^{l}_{g_{k+1}}(g_{k+1} - \phi^\ast g)|_{g_{k+1}} = O(r^{-3-l})
\]
for every $l\geq 0$.
\item We say that $M$ is of type $D_m$ for some $m\geq 0$ if there exists a compact set $K\subset M$, $R>0$ and a double cover $\phi\co H^{2m-4} \ra M\setminus K$ such that the group $\Z_2$ of deck transformations is generated by the standard involution on $H^{2m-4}$ and
\[
|\nabla ^{l}_{g_{2m-4}}(g_{2m-4} - \phi^\ast g)|_{g_{2m-4}} = O(r^{-3-l})
\]
for every $l\geq 0$.
\end{enumerate}
\end{definition}

By \cite[Theorem 1.1]{Chen:Chen:II} every ALF gravitational instanton is either of type $A_k$ for some $k\geq -1$ or $D_m$ for some $m\geq 0$ (the constraints $k\geq -1$ and $m \geq 0$ follow from \cite[Theorem 0.1]{Minerbe:Mass} and \cite[Corollary 3.2]{Biquard:Minerbe}, respectively).

\begin{remark*}
In the cyclic case Chen--Chen \cite[Theorem 1.1]{Chen:Chen:II} have a worse decay $O(r^{-2})$ of an ALF metric of type $A_k$ to $g_{k+1}$. However, from the explicit description of cyclic ALF gravitational instantons as multi-Taub--NUT spaces, as we will recall below, it is clear that one can always change coordinates by a translation on $\R^3$ so that every $A_k$ ALF space satisfies the stronger decay stated in Definition \ref{def:Ak:Dk}.(i).
\end{remark*}

\subsubsection{ALF spaces of cyclic type}

We saw that gravitational instantons of type $A_k$ can be constructed from Dirac monopoles on $\R^3$ with $k+1$ singularities via the Gibbons--Hawking ansatz. These are usually called \emph{multi-Taub--NUT} metrics. The case $k=0$ is the Taub--NUT metric on $\R^4$ and $k=-1$ is $\R^3\times \Sph^1$ with its flat metric. Minerbe \cite[Theorem 0.2]{Minerbe:Ak} has shown that every ALF space of cyclic type must be isometric to a multi-Taub--NUT metric.

From their explicit description one can easily compute basic information about cyclic ALF spaces: the fundamental group $\pi_1 (M)$, the second Betti number $b_2 (M)$, the Euler characteristic and the dimension of the moduli space $\mathcal{M}$ of $A_k$ metrics:
\begin{center}
\begin{tabular}{|c|c|c|c|c|}
  \hline
$k$  & $\pi_1(M)$ & $b_2 (M)$ & $\chi (M)$ & $\text{dim}(\mathcal{M})$ \\ \hline 
$-1$ & $\Z$ & $0$ & $0$ & $0$\\ \hline

$k>-1$ & $1$ & $k$ & $k+1$ & $3k$\\
  \hline
\end{tabular}
\end{center}
Here we assume that the asymptotic length of the circle fibre is normalised to be $1$ so that $\text{dim}(\mathcal{M})$ does not include rescalings.

\subsubsection{ALF spaces of dihedral type}

ALF metrics of dihedral type are not globally given by the Gibbons--Hawking construction and in most cases are not explicit. A number of different constructions have appeared over the past 30 years, but only recently Chen--Chen \cite[Theorem 1.2]{Chen:Chen:II} have shown that all these constructions yield equivalent families of ALF metrics. We distinguish the cases $m=0,1,2$ and $m \geq 3$.

\begin{itemize}
\item[$m=0$:] The $D_0$ ALF manifold is the moduli space of centred charge $2$ monopoles on $\R^3$ with its natural $L^2$--metric, known as the \emph{Atiyah--Hitchin manifold}. The metric admits a cohomogeneity one isometric action of $SU(2)$ and is explicitly given in terms of elliptic integrals \cite[Chapter 11]{Atiyah:Hitchin}. The Atiyah--Hitchin manifold is diffeomorphic to the complement of a Veronese $\R\PP^2$ in $\Sph^4$ and therefore it retracts to $\R\PP^2$. The Atiyah--Hitchin metric does not admit deformations as a $D_0$ ALF metric except for scaling.
\item[$m=1$:] The double cover of the Atiyah--Hitchin manifold is a $D_1$ ALF space. As a smooth manifold it is diffeomorphic to the complement of $\R\PP^2$ in $\C\PP^2$, or equivalently to the total space of $\mathcal{O}(-4)$ over $S^2$. Exploiting the rotational invariance of the metric it can be shown \cite[Proposition 5.5]{Micallef:Wolfson:1} that the $2$--sphere in the interior is a strictly stable minimal sphere which is not holomorphic with respect to any complex structure compatible with the metric, a fact that will play a crucial role in the proof of Theorem \ref{thm:Minimal:Surfaces}. This rotationally invariant $D_1$ ALF metric admits a $3$--dimensional family of $D_1$ ALF deformations, sometimes referred to as the \emph{Dancer metrics}.
\end{itemize}

\begin{remark*}
The fact that the double cover of the Atiyah--Hitchin manifold admits a $3$--parameter family of $D_1$ ALF deformations can also be shown using methods similar to the ones developed in this paper. Indeed, it is known \cite[\S 5.4]{Hitchin:L2:cohomology} that the rotationally invariant $D_1$ ALF metric admits a unique $L^2$--integrable (in fact, exponentially decaying) anti-self-dual harmonic form $\eta$. This form yields a $3$--dimensional space of infinitesimal hyperk\"ahler deformations and an extension of the analysis needed for the proof of Theorem \ref{thm:Main:precise} could be used to integrate these infinitesimal deformations to genuine $D_1$ ALF metrics. In fact Dancer \cite{Dancer} has constructed a $3$--parameter family of hyperk\"ahler deformations of the rotationally invariant $D_1$ ALF metric using Nahm's equations and hyperk\"ahler quotient techniques: there exists a hyperk\"ahler $8$--manifold $\mathcal{N}$ constructed as a moduli space of solutions to Nahm's equations which admits a triholomorphic $U(1)$--action. Denote by $\mu\co \mathcal{N}\ra \R^3$ the corresponding hyperk\"ahler moment map. Dancer identifies the rotationally symmetric $D_1$ ALF metric with the hyperk\"ahler quotient $\mu^{-1}(0)/U(1)$. By varying the level set of the moment map he then obtains a $3$--parameter family of hyperk\"ahler deformations of the Atiyah--Hitchin metric. By a general formula for the infinitesimal deformation of the symplectic form of a symplectic quotient corresponding to varying the level set of the moment map \cite{Variation:Symplectic:Quotient}, the infinitesimal deformations of the Atiyah--Hitchin metric corresponding to Dancer's metrics coincide with those determined by the $L^2$ harmonic form $\eta$, which is interpreted in this context as the curvature of the natural hyperholomorphic connection on the $U(1)$--bundle $\mu^{-1}(0)\ra \mu^{-1}(0)/U(1)$ induced by the Levi--Civita connection of $\mathcal{N}$.
\end{remark*}

\begin{itemize}
\item[$m=2$:] $D_2$ ALF metrics were constructed by Hitchin \cite[\S 7]{Hitchin:Einstein} using twistor methods and by Biquard--Minerbe \cite[Theorem 2.4]{Biquard:Minerbe} using a non-compact version of the Kummer construction: one considers the quotient of $\R^3 \times \Sph^1$ by an involution and resolves the two singularities gluing in copies of the Eguchi--Hanson metric.
\end{itemize}

\begin{remark*}
Biquard--Minerbe \cite[Theorem 2.4]{Biquard:Minerbe} use singular perturbation methods to solve a complex Monge--Amp\`ere equation on the minimal resolution of $(\R^3\times \Sph^1)/\Z_2$. Using the more general approach adopted in this paper to glue hyperk\"ahler structures one could extend their construction to recover a $6$--dimensional family of $D_2$ ALF metrics. 
\end{remark*}

\begin{itemize} 
\item[$m\geq 3$:] $D_m$ ALF metrics (for all $m \geq 1$) appeared in the work of Cherkis--Kapustin \cite{Cherkis:Kapustin} on moduli spaces of singular monopoles on $\R^3$ and were rigorously constructed by Cherkis--Hitchin \cite{Cherkis:Hitchin} using twistor methods and the generalised Legendre transform. In the case $m \geq 3$ a more transparent construction due to Biquard--Minerbe \cite[Theorem 2.5]{Biquard:Minerbe} yields $D_m$ ALF metrics by desingularising the quotient of the Taub--NUT metric by the binary dihedral group $\mathcal{D}_m$ of order $4(m-2)$ using ALE dihedral spaces. Using complex Monge--Amp\`ere methods Auvray \cite{Auvray:I,Auvray:II} has then constructed $3m$--dimensional families of $D_m$ ALF metrics on the smooth $4$--manifold underlying the minimal resolution of $\C^2/\mathcal{D}_m$.
\end{itemize}

\begin{remark}\label{rmk:Dk}
The gluing construction presented in this paper could be extended to the non-compact setting to yield yet another construction of dihedral ALF metrics. Indeed, one considers a Gibbons--Hawking metric obtained from the harmonic function
\[
h=\lambda-\frac{2}{|x|}+\sum_{i=1}^{m}{\frac{1}{2|x-x_i|}+\frac{1}{2|x+x_i|}}
\]
for $m$ distinct points $x_1,\dots, x_m \in \R^3\setminus\{ 0\}$. Observe that for $\lambda>0$ sufficiently large $h>0$ outside an arbitrarily small neighbourhood of the origin. Since the configuration of punctures is invariant under the standard involution of $\R^3$, this (incomplete) metric descends to a hyperk\"ahler metric on a $\Z_2$ quotient. For $\lambda$ sufficiently large one can then complete this metric by gluing in a copy of the $D_0$ ALF space close to the origin. This approximate solution could then be deformed to an exact hyperk\"ahler metric in a way similar to the proof of Theorem \ref{thm:Main:precise}.
\end{remark}

We summarise some of the properties of a $D_m$ ALF gravitational instanton $M$ in the following table:
\begin{center}
\begin{tabular}{|c|c|c|c|c|}
  \hline
$m$  & $\pi_1(M)$ & $b_2 (M)$ & $\chi (M)$ & $\text{dim}(\mathcal{M})$ \\ \hline 
$0$ & $\Z_2$ & $0$ & $1$ & $0$\\ \hline

$m>0$ & $1$ & $m$ & $m+1$ & $3m$\\
\hline
\end{tabular}
\end{center}
As in the cyclic case, $\mathcal{M}$ is the moduli space of $D_m$ ALF metrics modulo scaling.
 
\section{Gibbons--Hawking ansatz on a punctured $3$--torus}\label{sec:GH:3:torus}

In this section we use the Gibbons--Hawking ansatz \eqref{eq:Gibbons:Hawking} to construct families $\{ g^{\textup{gh}}_\epsilon\} _{\epsilon>0}$ of (incomplete) hyperk\"ahler metrics on circle bundles over a punctured $3$--torus modulo an involution. The parameter $\epsilon>0$ essentially determines the length of the circle fibres. As $\epsilon\ra 0$ the metric $g^{\textup{gh}}_\epsilon$ collapses to the flat orbifold $T^3/\Z_2$ with bounded curvature away from the punctures. The family $\{ g^{\textup{gh}}_\epsilon\} _{\epsilon>0}$ will serve as the model for a family of hyperk\"ahler metrics on the K3 surface collapsing to a $3$--dimensional limit in the region where the collapsing occurs with bounded curvature. In the next section we will use ALF gravitational instantons as models for high curvature regions to extend the metric $g^\textup{gh}_\epsilon$ to a complete almost hyperk\"ahler metric.

\subsection{Dirac monopoles on a punctured torus}

Let $\T=\R^3/\Lambda$ be a $3$--torus for some lattice $\Lambda \simeq \Z^3$. Endow $\T$ with a flat metric $g_\T$.

Let $\tau\co \T\ra \T$ be the standard involution $x\mapsto -x$ on $\T$ and denote by $q_1,\dots, q_8$ its fixed points. For each $j=1,\dots, 8$ choose a non-negative integer $m_j$.

Fix a $\tau$--symmetric configuration of further $2n$ distinct points $p_1,\tau (p_1), \dots, p_n,\tau(p_n)$. Sometimes we will use the notation  $-p_i$ for $\tau(p_i)$. Denote by $\T^\ast$ the punctured torus
\[
\T^\ast =\T\setminus \{ q_1,\dots, q_8, p_1,\tau (p_1),\dots, p_n, \tau (p_n)\}.
\]
Finally choose integer weights $k_1, \dots, k_n >0$ and assume the following \emph{balancing condition} holds:
\begin{equation}\label{eq:Balancing}
\sum_{j=1}^8{m_j}+\sum_{i=1}^n{k_i}=16.
\end{equation}
In particular, $n \leq \sum_{i=1}^n{k_i} \leq 16$.

For each $j=1,\dots, 8$ let $\rho_j$ denote the distance function from the point $q_j$ with respect to $g_\T$. Similarly, by abuse of notation we let $\rho_i$ denote the distance function from $\pm p_i$ in $\T/\tau$. By restricting the branched double cover $\T \ra \T/\tau$ to a sufficiently small ball centred at $\pm p_i$ in $\T/\tau$ we will also regard $\rho_i$ as the distance function on $\T$ from the point $p_i$ or $\tau(p_i)$.  

We look for a Dirac monopole $(h,\theta)$ on $\T^\ast$ with the following singular behaviour: $h$ is a harmonic function on $\T^\ast$ with prescribed singularities at the punctures
\begin{equation}\label{eq:Dirac:singularities}
h \sim \frac{2m_j-4}{2\rho_j} \mbox{ as } \rho_j \ra 0, \qquad h \sim \frac{k_i}{2\rho_i} \mbox{ as } \rho_i \ra 0.
\end{equation} 

\begin{prop}\label{prop:Dirac:torus}
Assume the balancing condition \eqref{eq:Balancing} is satisfied.
\begin{enumerate}
\item There exists a harmonic function $h$ on $\T^\ast$ with prescribed singular behaviour \eqref{eq:Dirac:singularities} such that $\ast_{g_\T} dh$ is the curvature $d\theta$ of a connection $\theta$ on some principal $U(1)$--bundle $P\ra \T^\ast$.
\item The moduli space of Dirac monopoles $(h,\theta)$ on $P$ is isomorphic to $\R \times \hat{\T}$, where $\hat{\T}$ is the dual torus parametrising flat $U(1)$--connections on $\T$.
\item The involution $\tau$ lifts to the involution $\tilde{\tau}$ of the $U(1)$--bundle $P$ which acts simultaneously as $\tau$ on $\T^\ast$ and as the standard involution on the circle fibres.
\end{enumerate}
\proof
The necessary and sufficient condition for the existence of the harmonic function $h$ is
\begin{equation}
\sum_{j=1}^8{2m_j-4}+2\sum_{i=1}^n{k_i}=\frac{1}{2\pi}\int_{\partial \T_{\sigma}}{\ast dh}=0,
\end{equation}
where $\T_\sigma$ denotes the complement of the union of small balls of radius $\sigma$ centred at the punctures. Thus if \eqref{eq:Balancing} is satisfied, a harmonic function $h$ with the singular behaviour \eqref{eq:Dirac:singularities} does indeed exists and is unique up to the addition of a constant.

By Lefschetz--Poincar\'e duality $H_2(\T^\ast)\simeq H^1_{c}(\T^\ast)$. The latter group sits in a long exact sequence
\[
0 \ra H^0(\T)\ra \Z^{2n+8}\ra H^1_c(\T^\ast)\ra H^1(\T)\ra 0,
\]
where $\Z^{2n+8}$ is generated by the $2n+8$ punctures. Thus $H_2(\T^\ast)$ is $(2n+10)$--dimensional and maps onto $H_2(\T)$ with kernel spanned by the classes of $2n+8$ spheres centred at the punctures. Note that the sum of these $2n+8$ homology classes vanishes.

Because of \eqref{eq:Dirac:singularities}, $2n+7$ of the $2n+10$ integrality constraints on $\tfrac{i}{2\pi}\ast dh$ to represent the first Chern class of a line bundle are automatically satisfied since we chose $2m_j-4, k_i\in\Z$. The remaining $3$ constraints can be reinterpreted in terms of the position of the punctures following the arguments in the proof of \cite[Proposition 3.5]{Charbonneau:Hurtubise}:
\[
\sum_{j=1}^8{(2m_j-8)\, q_j} + \sum_{i=1}^n{k_i\, \big( p_i+\tau (p_i) \big)}\in\Lambda.
\]
Since the points $q_j$ belong to the half-lattice $\tfrac{1}{2}\Lambda$ this condition is automatically satisfied.

We have therefore proved the existence of a principal $U(1)$ bundle $P\ra \T^\ast$ endowed with a connection $\theta$ with curvature $\ast dh$. Since $\T$ is not simply connected $\theta$ is uniquely determined up to a flat connection, \ie a point of the dual torus $\hat{\T}$. This concludes the proof of (i) and (ii).

By uniqueness up to the addition of a constant the harmonic function $h$ is $\tau$--invariant and therefore can be thought of as defined on $\T^\ast/\tau$. Since $\tau^\ast (\ast dh) = -\ast dh$, we can lift $\tau$ (uniquely up to gauge transformations) to an involution $\tilde{\tau}$ of the circle bundle $P$ by requiring that $\tilde{\tau}$ acts simultaneously as $\tau$ on $\T^\ast$ and as the standard involution on the circle fibres.
\endproof
\end{prop}

\subsection{Collapsed $S^1$--invariant hyperk\"ahler metrics}

Fix once and for all a Dirac monopole $(h,\theta)$ amongst the ones produced by Proposition \ref{prop:Dirac:torus}. Via the Gibbons--Hawking ansatz \eqref{eq:Gibbons:Hawking} we now use $(h,\theta)$ to construct (incomplete) hyperk\"ahler metrics with a triholomorphic circle action with orbits of small length.

Fix a (small) positive number $\epsilon>0$ and define
\begin{equation}\label{eq:Harmonic}
h_\epsilon = 1+\epsilon h.
\end{equation}
The Gibbons--Hawking ansatz \eqref{eq:Gibbons:Hawking} yields a hyperk\"ahler structure on $P|_{\mathcal{U}_\epsilon}$, where $\mathcal{U}_\epsilon$ is the open set of $\T^\ast$ where $h_\epsilon>0$. The hyperk\"ahler triple $\triple{\omega}^{\textup{gh}}_\epsilon$ and the induced metric $g^\textup{gh}_\epsilon$ are
\begin{equation}\label{eq:GH:collapsed}
\omega^{\textup{gh}}_{\epsilon,i} = \epsilon\, \theta_i \wedge \theta + h_\epsilon\, \theta_j\wedge\theta_k, \qquad g^{\textup{gh}}_\epsilon = h_\epsilon\,  \pi^\ast g_\T + \epsilon^2 h_\epsilon^{-1}\,\theta^2. 
\end{equation}
Here $(\theta_1,\theta_2,\theta_3)$ is a triple of closed $1$--forms on $\T$ such that $g_\T = \theta_1 ^2 + \theta_2 ^2 + \theta_3 ^2$. Note that the hyperk\"ahler structure $\triple{\omega}^{\textup{gh}}_\epsilon$ is $\tilde{\tau}$--invariant and therefore defines an induced hyperk\"ahler structure on the quotient $M^{\textup{gh}}_\epsilon = \left( P|_{\mathcal{U}_\epsilon} \right)/\tilde{\tau}$. For ease of notation, we will denote the induced hyperk\"ahler triple and metric on $M^\textup{gh}_\epsilon$ with the same symbols.

In the rest of the section we study the properties of the hyperk\"ahler manifold $M^{\textup{gh}}_\epsilon$. We aim to (i) study the local structure of the metric $g^{\textup{gh}}_\epsilon$ close to the punctures, (ii) determine the set where $h_\epsilon>0$, and (iii) understand the limit of $g^{gh}_\epsilon$ as $\epsilon\ra 0$. 

The following asymptotic expansions for the harmonic function $h_\epsilon$ close to the punctures are standard. In the case of a fixed point of the involution $\tau$ the improved decay follows from the fact that linear harmonic functions on $\R^3$ are not $\Z_2$--invariant.

\begin{lemma}\label{lem:harmonic:behaviour:punctures}
There exists $0<\rho_0<\tfrac{1}{4}\textrm{inj}\,g_\T$ such that the balls $B_{2\rho_0}(q_j)$, $j=1,\dots, 8$, and $B_{2\rho_0}(\pm p_i)$, $i=1,\dots, n$, in $\T$ are all disjoint and such that the following holds.
\begin{enumerate}
\item For $j=1,\dots,8$ there exists $\lambda_j\in\R$ such that in $B_{2\rho_0}(q_j)$
\[
h_\epsilon = (1+\epsilon\lambda_j) + \frac{\epsilon(m_j -2)}{\rho_j} + O(\epsilon\, \rho_j^2).
\]
\item For each $i=1,\dots,n$ there exists $\lambda_i\in\R$ and a linear function $\ell_i$ on $\R^3$ with $|\ell_i| \leq C\rho$ such that in $B_{2\rho_0}(\pm p_i)$
\[
h_\epsilon = (1+\epsilon\lambda_i) + \frac{\epsilon k_i}{2\rho_i} + \epsilon\,\ell_i + O(\epsilon\, \rho_i^2).
\]
\end{enumerate}
Moreover, $\rho_0, \lambda_j, \lambda_i, \ell_i$ depend continuously on the position of the punctures $p_1,\dots, p_n$ and on the flat metric $g_\T$ and $f=O(\epsilon\, \rho^3)$ means that there exists a constant $C$ depending continuously on these data such that $|\nabla ^k f|\leq C\epsilon \rho^{3-k}$ for $k=0,1,2,3$.
\end{lemma}

Since we chose $k_i>0$ for $i=1,\dots,n$ certainly a punctured neighbourhood of $\pm p_i$ is contained in the set $\mathcal{U}_\epsilon$ where $h_\epsilon>0$. As already mentioned, the Gibbons--Hawking metric can be extended by adding a single point to a smooth \emph{orbifold} metric modelled on $\C^2/\Z_{k_i}$. The obvious way to smooth out such an orbifold singularity is to replace the ``multiplicity'' $k_i$ point $p_i$ with $k_i$ points each with weight $1$. Then the Gibbons--Hawking ansatz yields a smooth metric that is modelled on a rescaled Taub--NUT space in a neighbourhood of $p_i$. However we prefer to leave the freedom to choose $k_i>1$ so that we can consider configurations of punctures that ``degenerate'' as $\epsilon \ra 0$ and see an $A_{k_i-1}$ ALF space appearing as a rescaled limit.

\begin{remark*}
One could also consider (but we will not pursue this in the paper) more general degenerating families of punctures with various clusters of points coalescing at different rates as $\epsilon \ra 0$. One would expect ``bubble trees'' of ALF \emph{and} ALE spaces appearing as rescaled limits in this case, \cf \cite[Remark 5.2]{Anderson:L2:structure}.
\end{remark*}

Next, we consider the structure of $g^{\textup{gh}}_\epsilon$ in a neighbourhood of $q_j$. By Lemma \ref{lem:harmonic:behaviour:punctures}.(i) certainly $h_\epsilon$ is positive in a punctured neighbourhood of $q_j$ whenever $m_j>2$. In this case the Gibbons--Hawking metric on $M^{\textup{gh}}_\epsilon$ can be extended to a smooth orbifold metric with a singularity of the form $\C^2/\mathcal{D}_{m_j}$, where $\mathcal{D}_{m_j}$ is the binary dihedral group of order $4(m_j-2)$. In contrast with the previous case, there is no explicit way to remove this singularity, but for \emph{fixed} $\epsilon>0$ one can imagine using the methods of \cite[\S 2.4]{Biquard:Minerbe} to resolve this singularity by gluing in a rescaled ALE dihedral space. However, we are interested in the limit $\epsilon\ra 0$ and in the next section we will directly glue in a $D_{m_j}$ ALF space to resolve this singularity.

Similarly, when $m_j=2$ one can choose $\epsilon$ sufficiently small so that $1+\epsilon\lambda_j>0$. Note that in this case $h_\epsilon$ and the $U(1)$--bundle $P$ are well defined at $q_j$. After quotienting by $\tilde{\tau}$, the Gibbons--Hawking metric $g^{\textup{gh}}_\epsilon$ becomes an orbifold metric modelled on $(\R^3\times S^1)/\Z_2$. As before, the two orbifold singularities could be resolved by (i) fixing $\epsilon>0$ and gluing in two copies of the Eguchi--Hanson metric, or (ii) letting $\epsilon\ra 0$ and gluing in a single copy of a $D_2$ ALF metric. We will follow the second approach. 

\begin{remark}\label{rmk:Kummer}
Note that if $m_j=2$ for all $j=1,\dots,8$ then $n=0$ by \eqref{eq:Balancing} and the bundle $P$ extends over every puncture: in fact $P$ is a $4$--torus and our construction reduces to the usual Kummer construction along a family of $4$--tori collapsing to a $3$--dimensional torus. This is the case considered by Page in \cite{Page:D2:ALF}. 
\end{remark} 

It remains to study the case when $m_j=0,1$ for some $j$. Assume this is the case for $j=1,\dots, k$ for some $1\leq k \leq 8$. Since $\sum_{j=1}^8{m_j}=16-\sum_{i=1}^n{k_i}\leq 16-n$, note that $k\geq 1$ as soon as $n\geq 1$ or $n=0$ and $m_j \neq 2$ for some $j$, \ie in every case except for the usual Kummer construction. The case $m_j=0,1$ is ``bad'' in the sense that $h_\epsilon\ra -\infty$ as $\rho_j\ra 0$.

\begin{lemma}\label{lem:Sign:harmonic:fct}
There exists $\epsilon_0>0$ depending continuously on $p_1,\dots,p_n$ and $g_\T$ such that for every $\epsilon<\epsilon_0$ we have $h_\epsilon >\tfrac{1}{2}$ on the complement of $\bigcup_{j=1}^k{B_{8\epsilon}(q_j)}$.
\proof
Restrict attention to the ball $B_{2\rho_0}(q_j)$. First note that $1+\frac{\epsilon\, (m_j-2)}{\rho}\geq 1-\frac{2\epsilon}{\rho} =\tfrac{3}{4}$ for $\rho= 8\epsilon$. Now choose $\epsilon_0>0$ so that $\epsilon (\lambda_j + C\tfrac{\epsilon^2}{64})\leq \tfrac{1}{4}$ for all $\epsilon\leq \epsilon_0$. Here $\lambda_j,C$ are the constants of Lemma \ref{lem:harmonic:behaviour:punctures}.(i). We conclude that $h_\epsilon>\tfrac{1}{2}$ on $\partial B_{8\epsilon}(q_j)$ for $\epsilon < \epsilon_0$. Since $h_\epsilon$ blows up to $+\infty$ at the punctures $q_{k+1}, \dots, q_8, \pm p_1,\dots, \pm p_n$ the maximum principle completes the proof.
\endproof
\end{lemma}

In other words, by choosing $\epsilon$ small enough we can assume that $h_\epsilon>0$ outside an arbitrarily small neighbourhood of the points $q_1,\dots, q_k$ where $m_j=0,1$.

Finally, we consider the limit $\epsilon\ra 0$.

\begin{lemma}\label{lem:GH:collapse}
As $\epsilon\ra 0$ the harmonic function $h_\epsilon$ converges to the constant function $1$. The convergence is in $C^{k,\alpha}$ on the complement of the union of balls of radius $\epsilon^{\beta}$ around the punctures, where $\beta>0$ is any number such that $\beta ^{-1}>k+1+\alpha$. In particular, $M^\textup{gh}_\epsilon$ collapses to the flat orbifold $\T/\tau$ with bounded curvature away from the punctures.
\proof
The first statement is a simple application of Lemma \ref{lem:harmonic:behaviour:punctures}, since close to each puncture we have
\[
\rho^k|\nabla ^k (h_\epsilon-1)|\leq C\epsilon\rho^{-1},
\]
where $\rho$ is the distance from the puncture. It follows that away from the punctures $g^{\textup{gh}}_\epsilon$ is $C^{k,\alpha}_{loc}$--close to the $\Z_2$ quotient of $g_\T + \epsilon^2\theta^2$ for any $k\geq 0$ and $\epsilon$ sufficiently small.
\endproof
\end{lemma}

\section{Approximate hyperk\"ahler metrics}\label{sec:Approximate:solution}

In this section we patch together ALF gravitational instantons and the incomplete hyperk\"ahler structure $\triple{\omega}^\textup{gh}_\epsilon$ of the previous section to construct a closed definite triple $\triple{\omega}_\epsilon$ which is approximately hyperk\"ahler. In the next section we will use analysis to deform $\triple{\omega}_\epsilon$ into a genuine hyperk\"ahler structure for each $\epsilon>0$ sufficiently small.

\subsection{The $4$--manifold $M_\epsilon$}

Let $\triple{\omega}^\textup{gh}_\epsilon$ be the hyperk\"ahler triple defined in \eqref{eq:GH:collapsed}. By Lemma \ref{lem:Sign:harmonic:fct} for $\epsilon>0$ small enough we think of $\triple{\omega}^\textup{gh}_\epsilon$ as defined on $M^\textup{gh}_\epsilon$, a smooth manifold with boundary obtained by restricting the line bundle $P$ to the complement of (arbitrarily) small balls centred at the punctures and then taking the quotient by the involution $\tilde{\tau}$. The boundary of $M^\textup{gh}_\epsilon$ has $n+8$ components, each of which has a collar neighbourhood diffeomorphic either to $H^{2m_j-4}/\Z_2$ for $j=1,\dots,8$, or $H^{k_i}$ for $i=1,\dots,n$.   

For each $j=1,\dots,8$ let $M_j$ be the smooth $4$--manifold underlying a $D_{m_j}$ ALF space and for each $i=1,\dots,n$ let $N_i$ be the smooth manifold underlying an $A_{k_i-2}$ ALF space. We construct a smooth $4$--manifold $M_\epsilon$ by cutting the ends of $M_j$ and $N_i$ and gluing the resulting manifolds with boundary to $M^\textup{gh}_\epsilon$ in a neighbourhood of $q_j$ or $\pm p_i$, respectively.

We will construct an approximate hyperk\"ahler structure on $M_\epsilon$ in the next subsection. Here we pause for a moment to determine the Betti numbers of $M_\epsilon$. While we will not use this result in an essential way in the rest of the paper, it is interesting to note how the balancing condition \eqref{eq:Balancing} appears naturally in the calculation of the Euler characteristic of $M_\epsilon$.

\begin{prop}\label{prop:Topology}
The Betti numbers of the compact orientable $4$--manifold $M_\epsilon$ are
\[
b_1 (M_\epsilon)=0, \qquad b_2^+(M_\epsilon)=3, \qquad b_2^- (M_\epsilon)=22.
\]
\proof
Decompose $M_\epsilon$ into the union of a piece $P/\tilde{\tau}$, an $A_{k_i-1}$ ALF space for each $i=1,\dots,n$ and a $D_{m_j}$ ALF space for each $j=1,\dots,8$. These pieces are identified along their common boundaries, which are homology spheres. Since all components have vanishing first Betti number, the reduced Mayer--Vietoris sequence yields $b_1 (M_\epsilon)=0$. The Euler characteristic is also easily calculated:
\[
\chi (M_\epsilon) = \chi (P/\Z_2) + \sum_{i=1}^n{\chi(A_{k_i-1})}+\sum_{j=1}^8{\chi (D_{m_j})} = 0+\sum_{i=1}^n{k_i}+\sum_{j=1}^8{m_j}+8=24
\]
by the balancing condition \eqref{eq:Balancing}.

It remains to calculate the signature $\tau(M_\epsilon)$. Below we will construct a definite triple $\triple{\omega}_\epsilon$ on $M_\epsilon$ which is close to define a hyperk\"ahler structure. By changing basis of $\Lambda^+ T^\ast M_\epsilon$ one can always deform this triple to a genuine $SU(2)$--structure (without requiring any differential constraint). In particular, $M_\epsilon$ can be endowed with an almost complex structure $J$ with $c_1(M_\epsilon,J)=0$. Since $c_2(M_\epsilon,J)=\chi(M_\epsilon)=24$, Hirzebruch's Signature Theorem and the equality of characteristic classes $p_1 =c_1^2-2c_2$ yield $\tau (M_\epsilon)=-16$.
\endproof
\end{prop}

\begin{remark*}
In Remark \ref{rmk:Kummer} we noted that the case where $n=0$ and $m_j=2$ for all $j=1,\dots,8$ reduces to the usual Kummer construction. Hence we know that $M_\epsilon$ is \emph{diffeomorphic} to the K3 surface in this special case. It seems likely one can prove that the diffeomorphism type of $M_\epsilon$ does not depend on the configuration of punctures satisfying the balancing condition \eqref{eq:Balancing}. Since we are going to construct a hyperk\"ahler metric on $M_\epsilon$, the calculation of the Betti numbers will anyway imply that $M_\epsilon$ is always diffeomorphic to the K3 surface.
\end{remark*}

\subsection{The definite triple $\underline{\omega}_\epsilon$}\label{sec:Perturbation:Ak}

We are now going to define an approximately hyperk\"ahler triple $\triple{\omega}_\epsilon$ on the $4$--manifold $M_\epsilon$.

For each $j=1,\dots,8$ denote by $\triple{\omega}_{q_j,\epsilon}$ the hyperk\"ahler triple obtained from the Gibbons--Hawking ansatz \eqref{eq:Gibbons:Hawking} using the harmonic function
\begin{equation}\label{eq:Harmonic:qj}
h_{q_j}=(1+\epsilon \lambda_j)+\frac{\epsilon(m_j-2)}{\rho}.
\end{equation}
By abuse of notation we think of $\triple{\omega}_{q_j,\epsilon}$ as defined both on the circle bundle $H^{2m_j-4}\ra \R^3 \setminus B_{8\epsilon}(0)$ as well as on its quotient by the involution that acts as the simultaneous standard involution on $\R^3$ and the circle fibres.

Similarly, for each $i=1,\dots,n$ let $\triple{\omega}_{p_i,\epsilon}$ be the hyperk\"ahler triple obtained from the Gibbons--Hawking ansatz using the harmonic function
\begin{equation}\label{eq:Harmonic:pi}
h_{p_i}=(1+\epsilon\lambda_i)+\frac{\epsilon\, k_i}{\rho} + \epsilon\, \ell_i.
\end{equation}
Here $\lambda_j,\lambda_i$ and $\ell_i$ are the constants and linear functions appearing in Lemma \ref{lem:harmonic:behaviour:punctures}. We will assume that $\epsilon$ is small enough to guarantee that $\tfrac{1}{2}<\lambda_i,\lambda_j <\tfrac{3}{2}$.

For all $j=1,\dots,8$ let $(M_j,\triple{\omega}_{M_j})$ be a complete $D_{m_j}$ ALF space. By Definition \ref{def:Ak:Dk} there exists a compact set $K\subset M_j$, $R_0>0$ and a diffeomorphism $M_j\setminus K \simeq H^{2m_j-4}/\Z_2$ such that
\[
\epsilon^2\triple{\omega}_{M_j} =\triple{\omega}_{q_j,\epsilon} + \triple{\eta}_{q_j,\epsilon}
\]
for $\rho>\epsilon R_0$ with $\triple{\eta}_{q_j,\epsilon}=O(\epsilon^3\rho^{-3})$ and similar estimates on the derivatives.

For each $i=1,\dots, n$ let $(N_i,\triple{\omega}_{N_i})$ be a complete $A_{k_i-1}$ ALF space. By the classification of ALF spaces of cyclic type \cite{Minerbe:Ak} the hyperk\"ahler structure on $N_i$ is explicitly given via the Gibbons--Hawking ansatz starting from a harmonic function on $\R^3$ with $k_i$ singularities. We can add to this function the smooth harmonic function $\epsilon^2 \ell_i$. Over a ball in $\R^3$ of radius much smaller than $\epsilon^{-2}$ we can regard the resulting hyperk\"ahler structure as a small perturbation of the ALF $A_{k_i-1}$ hyperk\"ahler structure $\triple{\omega}_{N_i}$. By abuse of notation we denote this perturbed hyperk\"ahler structure with the same symbol $\triple{\omega}_{N_i}$. The advantage of this modification is that now $\epsilon^2\triple{\omega}_{N_i}$ approaches $\triple{\omega}_{p_i,\epsilon}$ with a smaller error: by Definition \ref{def:Ak:Dk} there exists a compact set $K\subset N_i$, $R_0>0$ and a diffeomorphism $N_i\setminus K \simeq H^{k_i}|_{\R^3\setminus B_{R_0}}$ such that
\[
\epsilon^2\triple{\omega}_{N_i} =\triple{\omega}_{p_i,\epsilon} + \triple{\eta}_{p_i,\epsilon}
\]
with (by scaling) $\triple{\eta}_{p_i,\epsilon}=O(\epsilon^3\rho^{-3})$ and similar estimates on the derivatives.

\begin{remark*}
The choice of perturbing the gravitational instanton $N_i$ of type $A_{k_i-1}$ by adding a linear function on $\R^3$ can be regarded as an intermediate choice between resolving the singularity of $M^\textup{gh}_\epsilon$ by assuming all punctures $p_i$ have weight $k_i=1$ and the direct gluing of $N_i$ to $M^\textup{gh}_\epsilon$. While not strictly necessary, the choice of perturbing $N_i$ by a linear function makes the exposition more uniform. In particular, the closed definite triple we will construct below fails to be hyperk\"ahler by the \emph{same} amount in a neighbourhood of $q_j$ and $\pm p_i$.  
\end{remark*}

We will assume that $\epsilon$ is chosen so small as to make sure that $\epsilon R_0 \ll \rho_0$, where $\rho_0>0$ was fixed in Lemma \ref{lem:harmonic:behaviour:punctures}. By choosing $R_0$ larger if necessary we can assume that the harmonic functions $h_{q_j}$ and $h_{p_i}$ in \eqref{eq:Harmonic:qj} and \eqref{eq:Harmonic:pi} are as close to constant functions as we please for $\epsilon R_0 \leq \rho \leq \rho_0$. Then the hyperk\"ahler triples $\triple{\omega}_{q_j,\epsilon}$ and $\triple{\omega}_{p_i,\epsilon}$ define metrics $g_{q_j,\epsilon}$ and $g_{p_i,\epsilon}$ which are uniformly equivalent to $g_\T + \epsilon^2\theta^2$ in the regions $\epsilon R_0 \leq \rho_j, \rho_i \leq 2\rho_0$. In the rest of the section all norms and covariant derivatives will be computed with respect to this metric.

A crucial observation is that we can take $\triple{\eta}_{p_i,\epsilon}$ and $\triple{\eta}_{q_j,\epsilon}$ to be exact.

\begin{lemma}\label{lem:Asymptotics:Ak:Dk}
\begin{enumerate}
\item For all $i=1,\dots,n$ there exists a triple $\triple{a}_{p_i,\epsilon}$ of $1$--forms on $H^{k_i}$ such that
\[
|\nabla ^k \triple{a}_{p_i,\epsilon}| \leq C\epsilon^3\left( \frac{1}{\rho}\right) ^{2+k}
\]
and $\epsilon^2\triple{\omega}_{N_i} = \triple{\omega}_{p_i,\epsilon}+d\triple{a}_{p_i,\epsilon}$.
\item For all $j=1,\dots,8$ there exists a triple $\triple{a}_{q_j,\epsilon}$ of $\Z_2$--invariant $1$--forms on $H^{2m_j-4}$ such that
\[
|\nabla ^k \triple{a}_{q_j,\epsilon}| \leq C\epsilon^3\left( \frac{1}{\rho}\right) ^{2+k}
\]
and $\epsilon^2\triple{\omega}_{M_j} = \triple{\omega}_{q_j,\epsilon}+d\triple{a}_{q_j,\epsilon}$. 
\end{enumerate}
\proof
The proof is identical in the two cases. Set $k=k_i$ in case (i) and $k=2m_j-4$ in case (ii). In case (ii) we work with $\Z_2$--invariant forms on the double cover $H^{2m_j-4}$.

By scaling we can assume that $\epsilon=1$. It is enough to prove that every closed $2$--form $\eta$ with $\eta = O(\rho^{-3})$ can be written as $\eta = da$ with $|\nabla ^k a|=O(\rho^{-2-k})$.

Since the restriction of $H^{k}$ to an exterior domain in $\R^3$ is diffeomorphic to $(R,\infty)\times \Sigma$ with $\Sigma$ an homology sphere, we can write $\eta =d\rho \wedge\alpha + \beta$ for some $\rho$--dependent $1$--form $\alpha$ and $2$--form $\beta$ on $\Sigma$ with $|\alpha|+|\beta|=O(\rho^{-3})$.

The condition $d\eta=0$ implies $\partial_\rho\beta -d_\Sigma\alpha = 0$. We then define $a=-\int_{\rho}^\infty{\alpha}$. The Lemma follows.
\endproof
\end{lemma}

By a similar radial integration, Lemma \ref{lem:harmonic:behaviour:punctures}.(i) implies that in the regions $\epsilon R_0 \leq \rho_j \leq 2\rho_0$ and $\epsilon R_0 \leq \rho_i \leq 2\rho_0$, respectively, we can write
\begin{subequations}\label{eq:Asymptotics:punctures:gh}
\begin{equation}
\triple{\omega}^{\textup{gh}}_\epsilon = \triple{\omega}_{q_j,\epsilon} + d\triple{a}^{\textup{gh}}_{q_j,\epsilon}, \qquad \triple{\omega}^{\textup{gh}}_\epsilon = \triple{\omega}_{p_i,\epsilon} + d\triple{a}^{\textup{gh}}_{p_i,\epsilon}
\end{equation}
for triples $\triple{a}^{\textup{gh}}_{q_j,\epsilon}$ and $\triple{a}^{\textup{gh}}_{p_i,\epsilon}$ of $1$--forms such that
\begin{equation}
|\nabla ^k \triple{a}^{\textup{gh}}_{q_j,\epsilon}| \leq C\epsilon\rho_j^{3-k}, \qquad |\nabla ^k \triple{a}^{\textup{gh}}_{p_i,\epsilon}| \leq C\epsilon\rho_i^{3-k}
\end{equation}
\end{subequations}
for $k=0,1,2,3$.

Now, let $\chi_{q_j}$ and $\chi_{p_i}$ be cut-off functions with the following properties:
\begin{equation}\label{eq:cut:off}
\begin{gathered}
\chi_{q_j}\equiv 1 \text{ for } \rho_j \leq \epsilon^{\frac{2}{5}}, \qquad \chi_{q_j}\equiv 0 \text{ for } \rho_j \geq 2\epsilon^{\frac{2}{5}}, \qquad |\nabla \chi_{q_j}| \leq C\rho_j^{-1},\\
\chi_{p_i}\equiv 1 \text{ for } \rho_i \leq \epsilon^{\frac{2}{5}}, \qquad \chi_{p_i}\equiv 0 \text{ for } \rho_i \geq 2\epsilon^{\frac{2}{5}}, \qquad |\nabla \chi_{p_i}| \leq C\rho_i^{-1}.\\
\end{gathered}
\end{equation}

We finally define a triple of \emph{closed} $2$--forms on $M_\epsilon$ by
\begin{equation}\label{eq:approximate:triple}
\triple{\omega}_\epsilon = 
\begin{cases}
\epsilon^2\triple{\omega}_{M_j} & \mbox{if } \rho_j \leq \epsilon^{\frac{2}{5}},\\
\triple{\omega}_{q_j,\epsilon} + d\left( \chi_{q_j}\, \triple{a}_{q_j,\epsilon} + (1-\chi_{q_j})\,\triple{a}^{\textup{gh}}_{q_j,\epsilon}\right) & \mbox{if } \epsilon^{\frac{2}{5}}\leq \rho_j \leq 2\epsilon^{\frac{2}{5}},\\
\epsilon^2\triple{\omega}_{N_i} & \mbox{if } \rho_i \leq \epsilon^{\frac{2}{5}},\\
\triple{\omega}_{p_i,\epsilon} + d\left( \chi_{p_i}\, \triple{a}_{p_i,\epsilon} + (1-\chi_{p_i})\,\triple{a}^{\textup{gh}}_{p_i,\epsilon}\right) & \mbox{if } \epsilon^{\frac{2}{5}}\leq \rho_i \leq 2\epsilon^{\frac{2}{5}},\\
\triple{\omega}^{\textup{gh}}_{\epsilon} & \mbox{if } \rho_j \geq 2\epsilon^{\frac{2}{5}} \mbox{ and } \rho_i \geq 2\epsilon^{\frac{2}{5}} \mbox{ for all } i,j.
\end{cases}
\end{equation}

\begin{remark*}
When $k_i=1$ there is no need to glue in an $A_{0}$ ALF space ($\R^4$ endowed with the Taub--NUT metric), since $\triple{\omega}^\textup{gh}_\epsilon$ already extends smoothly over $\pm p_i$. It is however useful (we will use this in setting up the analysis for the deformation problem) to think of a rescaled Taub--NUT space localised around $\pm p_i$.
\end{remark*}

\subsubsection{The error}

We conclude this section by quantifying the failure of $\triple{\omega}_\epsilon = (\omega^1_\epsilon,\omega^2_\epsilon,\omega^3_\epsilon)$ to define a hyperk\"ahler structure. Since by construction $d\omega^i_\epsilon=0$ for all $i=1,2,3$, we only have to check that $\triple{\omega}_\epsilon$ is a definite triple and estimate the difference between the associated intersection matrix and the identity.

In the regions $\rho_j \leq \epsilon^{\frac{2}{5}}$, $\rho_i \leq \epsilon^{\frac{2}{5}}$ and when $\rho_i, \rho_j \geq 2\epsilon^{\frac{2}{5}}$ for all $i,j$ the triple $\triple{\omega}_\epsilon$ defines a genuine hyperk\"ahler structure. In the \emph{transition regions} $\epsilon^{\frac{2}{5}}\leq \rho_j \leq 2\epsilon^{\frac{2}{5}}$ and $\epsilon^{\frac{2}{5}}\leq \rho_i \leq 2\epsilon^{\frac{2}{5}}$ we have, respectively,
\[
\triple{\omega}_\epsilon - \triple{\omega}_{q_j,\epsilon}=O(\epsilon^{2-\frac{1}{5}}), \qquad \triple{\omega}_\epsilon - \triple{\omega}_{p_i,\epsilon}=O(\epsilon^{2-\frac{1}{5}}).
\]
Since $\triple{\omega}_{q_j,\epsilon}$ and $\triple{\omega}_{p_i,\epsilon}$ are hyperk\"ahler triples, we conclude that $\triple{\omega}_\epsilon$ is a definite triple for $\epsilon$ sufficiently small.

Let $\mu_\epsilon$, $g_\epsilon$ and $Q_\epsilon$ be the volume form, metric and intersection matrix associated to the definite triple $\triple{\omega}_\epsilon$ as in Section \ref{sec:Definite:triples}. Using Lemma \ref{lem:Asymptotics:Ak:Dk}, \eqref{eq:Asymptotics:punctures:gh}, \eqref{eq:cut:off} and the definition of $\triple{\omega}_\epsilon$ we calculate that
\begin{equation}\label{eq:Error}
|Q_\epsilon-\text{id}| \leq C \epsilon^{2-\frac{1}{5}}
\end{equation}
in every transition region $\epsilon^{\frac{2}{5}}\leq \rho_j \leq 2\epsilon^{\frac{2}{5}}$, $j=1,\dots,8$, and $\epsilon^{\frac{2}{5}}\leq \rho_i \leq 2\epsilon^{\frac{2}{5}}$, $i=1,\dots,n$. Outside the transition regions $Q_\epsilon\equiv \text{id}$. 

\section{Perturbation to hyperk\"ahler metrics}\label{sec:Deformation}

In the previous section we have constructed a $4$--manifold $M_\epsilon$ together with a closed definite triple $\triple{\omega}_\epsilon$ which is approximately hyperk\"ahler in the sense that the intersection matrix $Q_\epsilon$ associated with $\triple{\omega}_\epsilon$ differs from the identity by arbitrarily small terms as $\epsilon\ra 0$. We would like to deform $\triple{\omega}_\epsilon$ into a genuine hyperk\"ahler triple using analysis. Since the geometry of $M_\epsilon$ degenerates as as $\epsilon\ra 0$ we need to take some care in applying the Implicit Function Theorem.

As explained in Section \ref{sec:Definite:triples} we can reformulate the problem in terms of an elliptic PDE. Let $g_\epsilon$ be the Riemannian metric on $M_\epsilon$ defined by $\triple{\omega}_\epsilon$. Denote by $\mathcal{H}^+_\epsilon$ the space of self-dual harmonic forms with respect to $g_\epsilon$. By Proposition \ref{prop:Topology} $\mathcal{H}^+_\epsilon$ is $3$--dimensional spanned by $\omega^i_\epsilon$, $i=1,2,3$. The equation we want to solve is \eqref{eq:Nonlinear}, \ie
\begin{equation}\label{eq:Nonlinear:epsilon}
d^+\triple{a}+\triple{\zeta}=\mathcal{F}\left( (\text{id}-Q_\epsilon)-d^-\triple{a}^-\ast d^-\triple{a}\right), \qquad d^\ast\triple{a}=0,
\end{equation}
for a triple $\triple{a}$ of $1$--forms on $M_\epsilon$ and a triple $\triple{\zeta} \in \mathcal{H}^+_\epsilon\otimes \R^3$. 

The linearisation of \eqref{eq:Nonlinear:epsilon} is an isomorphism since $b_1 (M_\epsilon)=0$ by Proposition \ref{prop:Topology}. Our main task is to show that its inverse has bounded norm as $\epsilon\ra 0$ and to control the non-linearities in \eqref{eq:Nonlinear:epsilon} by introducing appropriate Banach spaces.

\subsection{The linear operator $d^\ast + 2\, d^+$ for collapsing Gibbons--Hawking metrics}
  
Before introducing weighted H\"older spaces and proving the main estimates, it is helpful to look more closely at the linearisation of \eqref{eq:Nonlinear:epsilon}. It involves the operator $D = d^\ast + 2\, d^+\co \Omega^1 (M_\epsilon) \ra \Omega^0 (M_\epsilon)\oplus \Omega^+(M_\epsilon)$, where adjoints and projections are computed using the metric $g_\epsilon$.

We are interested in understanding the behaviour of the operator $D$ (in particular, the presence of small eigenvalues) for a sequence of collapsing metrics in the Gibbons--Hawking form \eqref{eq:Gibbons:Hawking:metric}.

Consider then the metric
\[
g^\textup{gh}_\epsilon = h_\epsilon\, g_\T + \epsilon^2 h_\epsilon^{-1}\theta^2
\]
of \eqref{eq:GH:collapsed} over the circle bundle $\pi\co P\ra \mathcal{U}_\epsilon\subset \T^\ast$. By Lemma \ref{lem:Sign:harmonic:fct} the open sets $\mathcal{U}_\epsilon$ form an exhaustion of $\T^\ast$ as $\epsilon\ra 0$.

We first consider the geometry of $g^\textup{gh}_\epsilon$ and in particular calculate its Levi--Civita connection.

As before let $\theta_1,\theta_2,\theta_3$ denote closed $1$--forms on $\T$ such that $g_\T = \theta_1^2 + \theta_2^2 +\theta_3^2$. Let $\xi_1,\xi_2,\xi_3$ be the dual vector fields with respect to $g_\T$. We will not distinguish between a vector tangent to $\T$ and its horizontal lift to $P$ with respect to the connection $\theta$. In particular, $[\xi_i,\xi_j]=-d\theta (\xi_i,\xi_j)\, \xi$ as vector fields on $P$. Finally, let $\xi$ be the vertical vector field normalised so that $\theta (\xi)=1$.

Since $\theta([\xi,\xi_i])=-d\theta (\xi,\xi_i)=0$ and $\pi_\ast [\xi,\xi_i]=[\pi_\ast \xi,\xi_i]=0$, we have $[\xi,\xi_i]=0$. The Koszul formula
\[
2\langle \nabla _X Y,Z\rangle = \langle [X,Y],Z\rangle -\langle [Y,Z],X\rangle + \langle [Z,X],Y\rangle + X\cdot \langle Y,Z\rangle  + Y\cdot \langle X,Z\rangle-Z\cdot \langle X,Y\rangle
\]
then allows to calculate the Levi--Civita connection $\nabla$ of the metric $g^\textup{gh}_\epsilon$:
\[
\begin{gathered}
\nabla _\xi \xi=-\tfrac{1}{4}\epsilon^2 \nabla^\T h_\epsilon^{-2}, \qquad \nabla _{\xi_i}\xi=-\tfrac{1}{2}h_\epsilon^{-1}(\xi_i\cdot h_\epsilon)\,\xi +\tfrac{1}{2}\epsilon^2 h_\epsilon^{-2}(\xi_i\lrcorner d\theta)^{\sharp_\T},\\
\nabla _\xi \xi_i=-\tfrac{1}{2}h_\epsilon^{-1}(\xi_i\cdot h_\epsilon)\,\xi +\tfrac{1}{2}\epsilon^2 h_\epsilon^{-2}(\xi_i\lrcorner d\theta)^{\sharp_\T}, \\
\nabla _{\xi_j} \xi_i = -\tfrac{1}{2}d\theta (\xi_i,\xi_j)\, \xi + \tfrac{1}{2}h_\epsilon^{-1}\left( (\xi_j\cdot h_\epsilon)\, \xi_i - (\xi_i\cdot h_\epsilon)\, \xi_j - \delta_{ij}\nabla^\T h_\epsilon\right).
\end{gathered}
\]
Here $\nabla^\T$ and ${}^{\sharp_\T}$ denote gradient and musical isomorphism with respect to the flat metric $g_\T$ and we used the fact that $h_\epsilon$ is $S^1$--invariant. Since $\nabla$ is a metric connection, we calculate the covariant derivatives of the $1$--forms $\theta, \theta_i$ by duality. Since we will need this later, we write out formulas for $\nabla_\xi \theta$ and $\nabla_\xi\theta_i$:
\begin{equation}\label{eq:LC:GH}
\nabla _\xi \theta = \tfrac{1}{2}h_\epsilon^{-1}dh_\epsilon, \qquad \nabla_\xi\theta_i = -\tfrac{1}{2}\epsilon^2 h_\epsilon^{-3}(\xi_i\cdot h_\epsilon)\, \theta + \tfrac{1}{2}\epsilon^2 h_\epsilon^{-2}(\xi_i\lrcorner d\theta).
\end{equation}

Now, the cotangent bundle of $M^\textup{gh}_\epsilon$ is trivial as it is spanned by $\theta,\theta_1,\theta_2,\theta_3$. Thus we can write every $1$--form $a$ as
\begin{equation}\label{eq:a:GH}
a=\epsilon\, a_0\, \theta + a_1\, \theta_1 + a_2\, \theta_2 + a_3\, \theta_3
\end{equation}
for functions $a_0,a_1,a_2,a_3$. Note that
\[
|a|^2_{g^\textup{gh}_\epsilon} = h_\epsilon\, |a_0|^2 + h_\epsilon^{-1}\left( |a_1|^2 + |a_2|^2 + |a_3|^2 \right).
\]

A direct computation using the fact that $\theta_i$ is closed for $i=1,2,3$ and $(h_\epsilon,\epsilon\, d\theta)$ is a solution of the monopole equation \eqref{eq:Monopole} with respect to the flat metric $g_\T$ shows that
\begin{equation}\label{eq:D:GH}
\begin{gathered}
d^\ast a = -h_\epsilon^{-1}\left( \sum_{i=1}^3{\xi_i\cdot a_i} + \tfrac{1}{\epsilon} h_\epsilon^2\, \xi\cdot a_0 \right),\\
2\, d^+a = \sum_{i=1}^3 {\left( \xi_i\cdot a_0 + h_\epsilon^{-1} (\xi_j\cdot a_k - \xi_k \cdot a_j) +h_\epsilon^{-1}(\xi_i\cdot h_\epsilon)\, a_0-\tfrac{1}{\epsilon}(\xi\cdot a_i)\right) \omega^\textup{gh}_{\epsilon,i}},
\end{gathered}
\end{equation}
where $\omega^\textup{gh}_{\epsilon,i}$, $i=1,2,3$, is the hyperk\"ahler triple $\triple{\omega}^{\textup{gh}}_\epsilon$ of \eqref{eq:GH:collapsed}.

By Fourier analysis along the circle fibres we define projections $\Pi_0$ and $\Pi_\perp$ onto $S^1$--invariant and oscillatory components of functions. Via the trivialisation \eqref{eq:a:GH} $\Pi_0$ and $\Pi_\perp$ extend to $1$--forms. Since $h_\epsilon$ is $S^1$--invariant we see from \eqref{eq:D:GH} that the operator $D$ respects this decomposition.

For any fixed $\tau\in (0,1)$ restrict attention to the region in $M^\textup{gh}_\epsilon$ where $\rho_j,\rho_i \geq c\,\epsilon^{\frac{1-\tau}{2}}$ for all $i=1,\dots, n$ and $j=1,\dots, 8$. Then
\[
\| h_\epsilon -1\| _{C^0} \leq C\epsilon^{\frac{1+\tau}{2}}, \qquad \| \nabla ^\T h_\epsilon \| _{C^0} \leq C\epsilon^\tau
\]
by Lemma \ref{lem:GH:collapse}. We conclude that the operator $D$ of \eqref{eq:D:GH} acting on $S^1$--invariant $1$--forms approaches the Dirac operator
\begin{equation}\label{eq:Dirac:torus}
D_0 \co \Omega^0 (\T)\oplus\Omega^1(\T) \ra \Omega^0 (\T)\oplus\Omega^1(\T), \qquad (f,\gamma) \mapsto (d^\ast \gamma, df+\ast d\gamma)
\end{equation}
of the flat torus $\T$.

Moreover, using the expressions \eqref{eq:LC:GH} for $\nabla_\xi \theta$ and $\nabla_\xi \theta_i$ we find
\[
\epsilon^2 h_\epsilon ^{-1}|\nabla a |^2_{g^\textup{gh}_\epsilon} \geq |\nabla _\xi a|^2_{g^\textup{gh}_\epsilon} \geq h_\epsilon |\xi\cdot a_0|^2 + h_\epsilon^{-1}\sum_{i=1}^3{|\xi\cdot a_i|^2} - Ch_\epsilon^{-4}|\nabla^\T h_\epsilon|^2_{g_\T} \left( h_\epsilon |a_0|^2 + h_\epsilon^{-1}\sum_{i=1}^3{| a_i|^2}\right).
\]
Thus in the region where $\rho_j,\rho_i \geq c\,\epsilon^{\frac{1-\tau}{2}}$ for some $\tau>0$ we have
\begin{equation}\label{eq:Oscillatory}
\| \Pi_\perp a \|_{C^{0,\alpha}(S^1)} \leq C\epsilon\, \| \nabla a \|_{C^{0,\alpha}(S^1)}
\end{equation}
on each fibre for all $\epsilon$ sufficiently small.

These two observations --- the convergence of the operator $D$ to the Dirac operator $D_0$ of the flat $3$--torus as $\epsilon\ra 0$ and the strong control of the oscillatory part of $1$--forms --- will be crucial in the rest of the section. We will exploit the same remarks when considering blow-downs of ALF gravitational instantons. Let $(M,g_M)$ be a complete ALF space. Given a sequence $R_i\ra \infty$ consider the blow down $R_i^{-2} g_M$. Since by Definition \ref{def:Ak:Dk} $R_i^{-2} g_M$ is asymptotic (up to a double cover in the dihedral case) to the Gibbons--Hawking metric
\[
\left( 1+R_i^{-1}\frac{k}{2\rho}\right)\, g_{\R^3} + R_i^{-2} \left( 1+R_i^{-1}\frac{k}{2\rho}\right)^{-1}\theta^2,
\]
as $i\ra \infty$ the behaviour of the operator $D$ with respect to the metric $R_i^{-2} g_M$ is the same as the one observed for the metric $g^\textup{gh}_\epsilon$ as $\epsilon\ra 0$ with the flat $\R^3$ in place of the flat $3$--torus $\T$. Namely, on the region $\rho \geq c\, R_i ^{-\frac{1-\tau}{2}}$ \eqref{eq:Oscillatory} holds with $\epsilon=R_i^{-1}$ and the operator $D$ converges to the Dirac operator $D_0$ of flat space $\R^3$. 

\begin{remark*}
The behaviour of natural differential operators (the Laplacian acting on $p$--forms, the Dirac operator) associated with Riemannian metrics collapsing with bounded curvature and diameter have been studied by many authors, \cf for example \cite{Lott,Lott:Dirac}. The concrete situation we are interested in is a simple case of this more general theory and it seemed more appropriate to exploit the explicit nature of the Gibbons--Hawking metric rather than appealing to these more general results. 
\end{remark*}

In order to control the growth of differential forms close to the punctures on $\T^\ast$ and on the end of an ALF space we will now introduce weighted H\"older spaces.
 
\subsection{Weighted H\"older spaces}

We work on the Riemannian $4$--manifold $(M_\epsilon,g_\epsilon)$ constructed in Section \ref{sec:Approximate:solution}. The aim of this subsection is to introduce weighted H\"older spaces and prove a weighted Schauder estimate for the operator $D=d^\ast + 2\, d^+$ associated with the metric $g_\epsilon$.

For $R_0$ sufficiently large and $\epsilon=\epsilon (R_0)$ sufficiently small we define a \emph{weight function} $\rho_\epsilon$ as follows: we set
\begin{equation}\label{eq:Weight}
\rho_\epsilon = \begin{cases}
 \epsilon & \mbox{if } \rho_j \leq R_0\epsilon,\\
 \rho_j & \mbox{if } 2R_0\epsilon \leq \rho_j \leq \rho_0,\\
 \epsilon & \mbox{if } \rho_i \leq R_0\epsilon,\\
 \rho_i & \mbox{if } 2R_0\epsilon \leq \rho_i \leq \rho_0,\\
 1 & \mbox{if } \rho_j, \rho_i \geq 2\rho_0 \mbox{ for all }j=1,\dots,8, i=1,\dots, n,
 \end{cases}
\end{equation}
and let $\rho_\epsilon$ interpolate smoothly and monotonically between the various regions. By abuse of notation we think of $\rho_\epsilon$ as defined both on $M_\epsilon$ and on $M^\textup{gh}_\epsilon$ or its double cover $P|_{\mathcal{U}_\epsilon}$.

\begin{definition}\label{def:Holder}
For each $\delta\in\R$, $k\in\Z_{\geq 0}$ and $\alpha\in (0,1)$ define the weighted H\"older norm  $C^{k,\alpha}_\delta$ by
\[
\| a \| _{C^{k,\alpha}_{\delta}} = \sum_{j=1}^k{ \| \rho_\epsilon^{-\delta+j}\nabla^j a\| _{C^0} } + \text{sup}_{d(x,y)<\text{inj}\, g
_\epsilon}{ \min{\left\{\rho_\epsilon(x)^{-\delta +k+\alpha},\rho_\epsilon(y)^{-\delta +k+\alpha}\right\}}\frac{|\nabla ^ka(x)-\nabla^k a(y)|}{|x-y|^\alpha} }.
\]
Here all norms and covariant derivatives are computed with respect to the metric $g_\epsilon$ and $\nabla^k a(x)$ and $\nabla^k a(y)$ are compared using parallel transport along the unique geodesic connecting $x$ and $y$. Similarly set $\| a \| _{C^0_\delta}=\| \rho^{-\delta}a\| _{C^0}$.
\end{definition}

The following simple estimate for products in $C^{0,\alpha}_{\delta-1}$ will be used to control the non-linearities.

\begin{lemma}\label{lem:Product}
For every $\delta<1$ there exists a constant $C>0$ independent of $\epsilon$ such that
\[
\| u\, v \| _{C^{0,\alpha}_{\delta-1}} \leq C\epsilon^{\delta-1}\| u \| _{C^{0,\alpha}_{\delta-1}} \| v \| _{C^{0,\alpha}_{\delta-1}}.
\]
\proof
From the definition of the $C^{0,\alpha}_{\delta-1}$--norm it is immediate to check that
\[
\| u\, v \| _{C^{0,\alpha}_{\delta-1}} \leq C\| \rho_\epsilon^{\delta-1}\| _{C^0}\| u \| _{C^{0,\alpha}_{\delta-1}} \| v \| _{C^{0,\alpha}_{\delta-1}}.
\]
Since $\rho_\epsilon \geq c\, \epsilon$ and $\delta-1<0$ the result follows.
\endproof
\end{lemma}

We now consider the operator $D=d^\ast + 2\,d^+$ acting on $1$--forms of class $C^{1,\alpha}_{\delta}$. We prove the following weighted Schauder estimate.

\begin{prop}\label{prop:Schauder}
For every $\delta\in\R$ there exists a constant $C>0$ independent of $\epsilon$ such that
\[
\| a \| _{C^{1,\alpha}_\delta} \leq C \left( \| Da\| _{C^{0,\alpha}_{\delta-1}} + \| a\| _{C^0_\delta}\right).
\]
\proof
In order to prove this estimate it is convenient to cut the manifold $M_\epsilon$ in various pieces and analyse the geometry separately in each of them. The global estimate follows by combining the ``local'' estimates obtained in each of these pieces.

Consider first the region $\rho_j \leq 2R_0\epsilon$ for some $j=1,\dots,8$. The rescaled metric $\epsilon^{-2}g_\epsilon$ is isometric to a compact region in the $D_{m_j}$ ALF space $(M_j,g_{M_j})$.

Now, given a $1$--form $a$ on $M_\epsilon$, restrict $a$ to the region $\rho_j \leq 2R_0\epsilon$ and define $\tilde{a}=\epsilon^{-1-\delta}a$. The standard Schauder estimate for the elliptic operator $D$ associated with the metric $g_{M_j}$ is
\[
\| \tilde{a} \|_{C^{1,\alpha}} \leq C \left( \|D\tilde{a}\| _{C^{0,\alpha}} + \| \tilde{a}\| _{C^0}\right).
\]
Since $|\tilde{a}|_{\epsilon^{-2}g_\epsilon}=\epsilon^{-\delta}|a|_{g_\epsilon}$ and the norms $\nabla\tilde{a}$ and $D\tilde{a}$ are related in a similar way to those of $\nabla a$ and $Da$, the weighted Schauder estimate follows immediately.

The same argument can be applied in the region $\rho_i \leq 2R_0\epsilon$: the role of $M_j$ is now played by a small perturbation (\cf the beginning of Section \ref{sec:Perturbation:Ak}) of the $A_{k_i-1}$ ALF space $N_i$. 

Consider now the transition region $R_0\, \epsilon \leq \rho_j \leq \rho_0$ for some $j=1,\dots, 8$. We can work on the double cover $H^{2m_j-4}$ and restrict to $\Z_2$--invariant forms. Scaling by $\epsilon$ as above we reduce to consider the restriction of $H^{2m_j-4}$ to the region $R_0 \leq \rho \leq \frac{\rho_0}{\epsilon}$ in $\R^3$ endowed with a metric
\[
g=\left( 1+\epsilon\lambda_j+\frac{m_j-2}{\rho} \right)g_{\R^3} + \left(1+\epsilon\lambda_j+\frac{m_j-2}{\rho}\right)^{-1}\theta^2 + O(\rho^{-3}) + O(\epsilon^3\rho^2).
\]
Moreover, after rescaling the weight function $\rho_\epsilon$ coincides with the radial function $\rho$ on $\R^3$.

Fix a number $\sigma\in (0,1)$. For each point $x$ let $\pi(x)$ be its image in $\R^3$ and set $R=\sigma \rho(x)$. Up to changing $R_0$ and $\rho_0$ into $(1-\sigma) R_0$ and $(1+\sigma) \rho_0$ we can assume that $B_R (\pi(x))$ is contained in the annulus $R_0 \leq \rho \leq \frac{\rho_0}{\epsilon}$ and that the restriction of $H^{2m_j-4}$ to this ball is trivial. We can then work on a ``square'' $B_R\times [-R,R]$ in the universal cover of $H^{2m_j-4}|_{B_R}$. Rescaling the metric by $R^{-2}$, applying standard Schauder estimates, rescaling back and multiplying by $R^{-\delta}$ we obtain
\[
\| a \| _{C^{1,\alpha}_\delta (B_R)} \leq C \left( \| Da\| _{C^{0,\alpha}_{\delta-1}(B_R) } + \| a \| _{C^0_\delta (B_R)}\right).
\]

The case of the region $R_0\epsilon\leq \rho_i\leq \rho_0$ is completely analogous.

Finally, in the region where $\rho_j,\rho_i\geq \tfrac{1}{2}\rho_0$ the weight function $\rho_\epsilon$ is uniformly equivalent to the constant $1$ and therefore weighted spaces coincide with standard H\"older spaces. Moreover the harmonic function $h_\epsilon$ is $C^{\infty}$--close to the constant $1$. The metric $g_\epsilon$ is therefore $C^{\infty}$--close to the metric $g_\infty = g_{\T} + \epsilon^2 \theta^2$. The Schauder estimate for forms supported in this region is immediate since we can restrict to small balls in the torus $\T$ on which the circle bundle $P$ is trivial and then work on the universal cover, which has bounded geometry.
\endproof
\end{prop}

\subsection{The linear estimate}

We can now prove the main result about the linearisation of \eqref{eq:Nonlinear:epsilon}: the operator $D$ has uniformly bounded inverse.

\begin{prop}\label{prop:Linearisation}
For $\epsilon$ sufficiently small and $\delta \in (-2,0)$ there exist $C$ independent of $\epsilon$ such that
\[
\| a \| _{C^{1,\alpha}_\delta} \leq C \| Da\| _{C^{0,\alpha}_{\delta-1}}.
\]
\proof
By contradiction assume that there exists a sequence $\epsilon_i\ra 0$ and $1$--forms $a_i$ on $M_{\epsilon_i}$ such that $\| a_i \| _{C^{1,\alpha}_\delta}=1$ but $\| Da_i \| _{C^{0,\alpha}_{\delta-1}}\ra 0$.

First of all we show that for every compact set $K$ in $M^{\textup{gh}}_\epsilon$ which does not contain any puncture we must have $\| a_i \| _{C^{1,\alpha}_\delta(K)}\ra 0$.

Over $K$ we can work on the double cover of $M^{\textup{gh}}_\epsilon$ and regard $a_i$ as $\Z_2$--invariant forms. Write $a_i = \epsilon_i f_i\, \theta + \gamma_i$, for a function $f_i$ and a $1$--form $\gamma_i$ such that $\xi\lrcorner \gamma_i=0$ (recall that $\xi$ is the vector field dual to $\theta$). Over $K$ we can also decompose $a_i = \Pi_0 a_i + \Pi_\perp a_i$. Observe that for any $0<\tau<1$ \eqref{eq:Oscillatory} implies that
\[
\epsilon_i^{\frac{1-\tau}{2}}\| \Pi_\perp a_i \| _{C^0_\delta}\leq \| \Pi_\perp a_i \| _{C^0_{\delta-1}} \leq C\epsilon_i \| a\| _{C^{1,\alpha}_{\delta-1}}
\]
provided $\rho_j,\rho_i \geq c\, \epsilon_i^{\frac{1-\tau}{2}}$. Since the gluing regions occur for $\rho_j\sim \epsilon_i^{\frac{2}{5}}$ and $\rho_i\sim \epsilon_i^{\frac{2}{5}}$ we can choose $\tau>0$ sufficiently small so that these assumptions are satisfied on $M^\textup{gh}_\epsilon$. Thus $\| \Pi_\perp a_i \| _{C^0_\delta(K)}\ra 0$. 

By the Arzel\'a--Ascoli Theorem we can therefore assume that $(f_i,\gamma_i)$ converges to $(f_0,\gamma)\in \Omega^0 (\T)\oplus \Omega^1(\T)$. By \eqref{eq:D:GH} $(f_0,\gamma)$ satisfies
\begin{equation}\label{eq:Linear:collapsed:limit}
\ast d\gamma+df_0=0=d^\ast \gamma
\end{equation}
on $\T$. The control on $\| a_i \| _{C^{1,\alpha}_\delta}$ guarantees that $|f_0|+|\gamma| \leq C \rho^{\delta}$ close to the punctures.

We want to conclude that $f_0$ is constant and $\gamma$ is a smooth harmonic $1$--form on $\T$. By trivialising the cotangent bundle of the $3$--torus $\T$ by harmonic $1$--forms $\theta_1,\theta_2,\theta_3$ we can write $\gamma=f_1 \theta_1 + f_2 \theta_2 + f_3 \theta_3$. Then $f_0,f_1,f_2,f_3$ are harmonic functions on the punctured $3$--torus with controlled blow-up rate at the punctures. Since $\delta >-2$, close to each puncture we must have $f_i=\lambda_i + c_i \rho^{-1} +O(\rho)$ for some constants $\lambda_i,c_i$. However, $c_i$ must vanish for all $i=0,1,2,3$ if $(f_0,f_i)$ is a solution of the first order system \eqref{eq:Linear:collapsed:limit} and not only of the second order PDE this implies. Hence the functions $f_i$ are bounded harmonic functions on $\T$ and must be constant.

However, by $\Z_2$--invariance of the $1$--forms $a_i$ the functions $f_i$ must be odd with respect to the involution $\tau$ on $\T$ and must therefore vanish. Thus $(f_0,\gamma)=0$ and the Schauder estimate of Proposition \ref{prop:Schauder} implies that $\| a_i \| _{C^{1,\alpha}_\delta (K)}\ra 0$.

Next we look at what happens close to one of the punctures. By what we have just proved and the assumption $\| a_i \|_{C^{1,\alpha}_\delta (M_\epsilon)}=1$, there exists at least a $j=1,\dots, 8$ or $i=1,\dots, n$ such that $\| a_i \|_{C^{1,\alpha}_\delta} \geq \tfrac{1}{n+8}>0$ in the region $\rho_j \leq \epsilon^\frac{2}{5}$ or $\rho_i \leq \epsilon^{\frac{2}{5}}$. We fix attention to such a region. Rescaling the metric by $\epsilon_i^{-2}$ and replacing $a_i$ with $\tilde{a}_i=\epsilon_i^{-\delta-1} a_i$, from now on we will work on a $D_m$ or an $A_{k-1}$ ALF space. Denote either of these non-compact manifolds by $M$. Because of the behaviour of the weighted H\"older norm in Definition \ref{def:Holder} under rescaling, we have $\| \tilde{a}_i \|_{C^{1,\alpha}_\delta } = \| a_i \|_{C^{1,\alpha}_\delta}\geq \tfrac{1}{n+8}$. For ease of notation we replace $\tilde{a}_i$ with $a_i$ until the end of the proof.

By the Arzel\'a--Ascoli Theorem over every compact set of $M$ we can extract a subsequence of $\{ a_i\}$ that converges to a solution $a$ of $Da=0$ on $M$. Moreover, $|a|\leq C\rho^\delta$. Since $\delta<0$ we conclude that $a=0$. Indeed, the equation $Da=0$ in particular implies that $\triangle a=0$. Since $M$ is Ricci-flat, the Weitzenb\"ock formula yields $|a|\,\triangle|a| \leq 0$ and therefore $|a|=0$ by the maximum principle.

Now, if $\| a_i \| _{C^0_\delta (M)}\ra 0$ then the Schauder estimate of Proposition \ref{prop:Schauder} would yield a contradiction to the assumption $\| a_i \| _{C^{1,\alpha}_\delta}\geq \frac{1}{n+8}$. Assume therefore that there exists some $\nu>0$ such that $\| a_i \| _{C^0_\delta}\geq \nu$. Since $a_i\ra 0$ in $C^{1,\alpha}_\delta(K)$ for every compact set $K\subset M$ there must exists a sequence of points $x_i\in M$ going off to infinity (in particular $R_i:=\rho(x_i)\ra \infty$) such that $|a_i (x_i)| \geq \nu \rho(x_i)^\delta$.

Now rescale the metric on $M$ by $R_i^{-2}$ and replace $a_i$ by $R_i^{-\delta-1}a_i$. Then $(M,R_i^{-2}g_M)$ is converging to the tangent cone $C$ at infinity of $M$, \ie either $C=\R^3$ or $C=\R^3/\Z_2$ depending on whether $M$ is of cyclic or dihedral type. 

As in the first step of the proof, we can use \eqref{eq:D:GH} and \eqref{eq:Oscillatory} to conclude that $R_i^{-\delta-1}a_i$ sub-converges over compact subsets of $C\setminus \{ 0\}$ to a pair $(f_0,\gamma)\in \Omega^0(C)\oplus \Omega^1 (C)$ such that
\[
\ast d\gamma+df_0=0=d^\ast \gamma,
\] 
$|f|+|\gamma|\leq C\rho^\delta$ and $(|f|+|\gamma|)(x_0)=\nu>0$ for some $x_0\in C\setminus \{ 0 \}$. As before, the fact that $\delta>-2$ implies that $f,\gamma$ are bounded close to the origin in $C$. The fact that $\delta<0$ then forces $(f_0,\gamma)$ to vanish. However this contradicts the fact that $(|f|+|\gamma|)(x_0)>0$.
\endproof
\end{prop}

\subsection{The non-linear problem}

We are now ready to deform the triple $\triple{\omega}_\epsilon$ into a genuine hyperk\"ahler triple by using the following Implicit Function Theorem.

\begin{lemma}\label{lem:IFT}
Let $\Phi\co E\ra F$ be the smooth function between Banach spaces and write $\Phi (x)=\Phi(0) + L(x)+N(x)$, where $L$ is linear and $N$ contains the non-linearities. Assume that there exists constants $r, C,q$ such that
\begin{enumerate}
\item $L$ is invertible with $\| L^{-1} \| \leq C$;
\item $\| N(x)-N(y)\|_F \leq q \|x+y\| _E \| x-y\| _E$ for all $x,y\in B_{r}(0)\subset E$;
\item $\| \Phi(0)\| _{F}< \min\left\{ \frac{r}{2C},\frac{1}{4qC^2}\right\}$.
\end{enumerate}
Then there exist a unique $x\in E$ with $\| x \| _E \leq 2C\| \Phi(0)\|_F$ such that $\Phi(x)=0$.
\end{lemma}

In our situation we set
\[
E:=\left( C^{1,\alpha}_\delta (T^\ast M_\epsilon) \oplus \mathcal{H}^+_\epsilon\right) \otimes \R^3,
\]
where $\mathcal{H}^+_\epsilon$ denotes the space of self-dual harmonic forms with respect to $g_\epsilon$, \ie constant linear combinations of $\omega^1_\epsilon, \omega^2_\epsilon,\omega^3_\epsilon$. We endow $E$ with the product of the $C^{1,\alpha}_\delta$--norm and the norm on the finite dimensional vector space $\mathcal{H}^+_\epsilon\otimes\R^3\simeq \R^9$ induced by the $L^2$--norm. Similarly we set
\[
F:=C^{0,\alpha}_{\delta-1} (\R \oplus \Lambda^+T^\ast M_\epsilon) \otimes \R^3
\]
endowed with the $C^{0,\alpha}_{\delta-1}$--norm.

The operator $\Phi$ is the one defined by \eqref{eq:Nonlinear:epsilon}. Thus $\Phi(0)=-\mathcal{F}(\text{id}-Q_\epsilon)$, $L(\triple{a}+\triple{\zeta})=D\triple{a}+\triple{\zeta}$ and the non-linear term is
\[
N(\triple{a}+\triple{\zeta})=\mathcal{F}\left( \text{id}-Q_\epsilon \right)-\mathcal{F}\left( \text{id}-Q_\epsilon - d^-\triple{a}\ast d^-\triple{a}\right).
\]
We need to check that the hypothesis of Lemma \ref{lem:IFT} are satisfied.

We use Proposition \ref{prop:Linearisation} to show that $L$ has uniformly bounded inverse for $\delta\in (-\tfrac{1}{2},0)$.

\begin{lemma}\label{lem:Linear}
For $\delta \in (-\tfrac{1}{2},0)$ and $\epsilon$ sufficiently small there exists a constant $C>0$ independent of $\epsilon$ such that for every triple of self-dual $2$--forms $\triple{\xi}\in C^{0,\alpha}_{\delta-1}$ there exists a unique $(\triple{a},\triple{\zeta})\in E$ with
\[
\| \triple{a}\| _{C^{1,\alpha}_\delta}+\| \triple{\zeta}\| \leq C \| \triple{\xi}\|_{C^{0,\alpha}_{\delta-1}}.
\]
and $L(\triple{a},\triple{\zeta})=\triple{\xi}$.
\proof
First of all, note that the $2$--forms $\omega_\epsilon^i$ have uniformly bounded $C^{0,\alpha}_{\delta-1}$--norm. Indeed, outside the gluing regions $\triple{\omega}_\epsilon$ is a hyperk\"ahler triple and thus $\omega_\epsilon^i$ is parallel and bounded. On the gluing regions, $\triple{\omega}_\epsilon$ differs from the hyperk\"ahler triple $\triple{\omega}_{q_j,\epsilon}$ or $\triple{\omega}_{p_i,\epsilon}$ by terms of order $O(\epsilon \rho^2 + \epsilon^3\rho^{-3})$ (with similar estimates on their derivatives). Finally, $\rho_\epsilon^{-\delta +1}$ is bounded above since $\delta<0$.

Now, let $\widetilde{\triple{\omega}}$ be an $L^2$--orthonormal triple of harmonic self-dual forms with respect to $g_\epsilon$. Since
\[
\int_{M_\epsilon}{\omega_\epsilon^i \wedge \omega_\epsilon^j} = 2\int_{M_\epsilon}{(Q_\epsilon)_{ij}\,\dvol_{g_\epsilon}},
\]
$Q_\epsilon$ is close to the identity and $\vol _{g_\epsilon} (M_\epsilon)=O(\epsilon)$ we can assume that
\[
\| \widetilde{\triple{\omega}}\| _{C^{0,\alpha}_{\delta-1}}\leq C \epsilon^{-\frac{1}{2}}\| \triple{\omega}_\epsilon \| _{C^{0,\alpha}_{\delta-1}}\leq C \epsilon ^{-\frac{1}{2}}.
\]

Finally, observe that for every $u\in C^{0,\alpha}_{\delta-1}$ we have
\[
\| u \| _{L^2} \leq \| \rho_\epsilon^{\delta -1}\|_{L^2}\| u \| _{C^{0,\alpha}_{\delta-1}} \leq C (\epsilon^{\frac{1}{2}} + \epsilon ^{\delta +1}) \| u \| _{C^{0,\alpha}_{\delta-1}}.
\]
Indeed, using the definition \eqref{eq:Weight} of $\rho_\epsilon$ and the construction of $\triple{\omega}_\epsilon$ it is not difficult to estimate $\| \rho_\epsilon^{\delta -1}\|_{L^2} \leq C (\epsilon^{\frac{1}{2}}+ \epsilon^{\delta +1})$.

Now let $\pi\co C^{0,\alpha}_{\delta-1} (\Lambda^+T^\ast M_\epsilon) \ra \mathcal{H}^+_\epsilon$ be the $L^2$--orthogonal projection
\[
\pi (\xi) = \sum_{i=1}^3{\lambda_i\, \widetilde{\omega}_i}, \qquad \lambda _i =\int{\xi \wedge \widetilde{\omega}_i},
\]
and regard $\text{id}-\pi$ as a map $C^{0,\alpha}_{\delta-1} (\Lambda^+T^\ast M_\epsilon) \ra C^{0,\alpha}_{\delta-1} (\Lambda^+T^\ast M_\epsilon)$. By the remarks above we have
\[
|\lambda_i|\leq C (1+\epsilon^{\delta +1})\| \xi \| _{C^{0,\alpha}_{\delta-1}}, \qquad \| \pi (\xi)\| _{C^{0,\alpha}_{\delta-1}} \leq C (1+\epsilon^{\delta + \frac{1}{2}})\| \xi \| _{C^{0,\alpha}_{\delta-1}}.
\]
Thus if $\delta \geq -\tfrac{1}{2}$ the projections $\pi$ and $\text{id}-\pi$ are uniformly bounded. Proposition \ref{prop:Linearisation} and the surjectivity of $L$ then yield the result. 
\endproof
\end{lemma}

Next, we consider the non-linear term $N$. Note that this does not involve the harmonic part $\triple{\zeta}$. The function $\mathcal{F}$ is pointwise smooth with uniformly controlled norm for $\epsilon$ sufficiently small. Using the Taylor expansion of $\mathcal{F}$ at $\text{id}-Q_\epsilon$ and Lemma \ref{lem:Product} to control products we can therefore find $r>0$ and $C$ independent of $\epsilon$ such that assumption (ii) in Lemma \ref{lem:IFT} is satisfied with $r$ and $q=C\epsilon^{\delta-1}$ for some $\epsilon$--independent constant $C$. 

Finally,
\[
\| \mathcal{F}(\text{id}-Q_\epsilon)\| _{C^{0,\alpha}_{\delta-1}} \leq C \| \text{id}-Q_\epsilon\| _{C^{0,\alpha}_{\delta-1}} \leq C \epsilon^{2+\frac{1}{5}-\frac{2}{5}\delta}.
\]
Indeed, setting $\rho=\rho_j$ for $j=1,\dots,8$ or $\rho =\rho_i$ for some $i=1,\dots,n$, in the region $\epsilon^{\frac{2}{5}} \leq \rho \leq 2\epsilon^{\frac{2}{5}}$ we have $|\text{id}-Q_\epsilon|=O(\epsilon \rho^2 + \epsilon^3\rho^{-3})$ by \eqref{eq:Error} and $\rho_\epsilon=\rho$ by \eqref{eq:Weight}.

Thus assumption (iii) in Lemma \ref{lem:IFT} is therefore satisfied as soon as $\epsilon^{2+\frac{1}{5}-\frac{2}{5}\delta}\ll \epsilon^{1-\delta}$, \ie
\[
\epsilon^{\frac{3}{5}(\delta+2)}\ll 1.
\]
If $\delta>-2$ this condition is satisfied for $\epsilon>0$ sufficiently small.

\begin{theorem}\label{thm:Main:precise}
Let $(\T,g_\T)$ be a flat $3$--torus with standard involution $\tau\co \T\ra \T$. Let $q_1,\dots, q_8$ be the fixed points of $\tau$ and let $p_1,\tau(p_1),\dots, p_n,\tau(p_n)$ be further $2n$ distinct points. Denote by $\T^\ast$ the punctured torus $\T\setminus \{ q_1,\dots,q_8,p_1,\dots,\tau(p_n)\}$.

Let $m_1,\dots, m_8\in\Z_{\geq 0}$ and $k_1,\dots,k_n\in\Z_{\geq 1}$ satisfy
\[
\sum_{j=1}^8{m_j}+\sum_{i=1}^n{k_i}=16.
\] 
For each $j=1,\dots,8$ fix a $D_{m_j}$ ALF space $M_j$ and for each $i=1,\dots,n$ an $A_{k_i-1}$ ALF space $N_i$.

Then there exists a $1$--parameter family of hyperk\"ahler metrics $\{ g_\epsilon\} _{\epsilon\in (0,\epsilon_0)}$ on the K3 surface with the following properties. We can decompose the K3 surface into the union of open sets $K^\epsilon \cup \bigcup_{j=1}^8{M_j^\epsilon}\cup\bigcup_{i=1}^n{N_i^\epsilon}$ such that
\begin{enumerate}
\item $(K^\epsilon,g_\epsilon)$ collapses to the flat orbifold $\T^\ast/\Z_2$ with bounded curvature away from the punctures;
\item for each $j=1,\dots, 8$ and $k\geq 0$, $(M_j^\epsilon,\epsilon^{-2}g_\epsilon)$ converges in $C^{k,\alpha}_{loc}$ to the $D_{m_j}$ ALF space $M_j$;
\item for each $i=1,\dots,n$ and $k\geq 0$, $(N_j^\epsilon,\epsilon^{-2}g_\epsilon)$ converges in $C^{k,\alpha}_{loc}$ to the $A_{k_i-1}$ ALF space $N_i$.
\end{enumerate}
\proof
Given data as in the statement we constructed a $4$--manifold $M_\epsilon$ and a $1$--parameter family of closed definite triples $\triple{\omega}_\epsilon$ which are approximately hyperk\"ahler. For $\epsilon$ sufficiently small we can apply Lemma \ref{lem:IFT} to find unique $\triple{a}_\epsilon\in C^{1,\alpha}_\delta (T^\ast M_\epsilon)$ for $\delta \in (-\tfrac{1}{2},0)$ and $\triple{\zeta}_\epsilon\in \mathcal{H}^+_\epsilon$ such that $\| \triple{a}_\epsilon \| _{C^{1,\alpha}_\delta}+\| \triple{\zeta}_\epsilon\| \leq C\epsilon^{\frac{11-2\delta}{5}}$ and $\triple{\omega}_\epsilon + d\triple{a}_\epsilon + \triple{\zeta}_\epsilon$ is a hyperk\"ahler structure on $M_\epsilon$. In particular, since $b_1 (M_\epsilon)=0$ by Proposition \ref{prop:Topology}, $M_\epsilon$ must be diffeomorphic to the K3 surface.

Away from the gluing regions $\triple{a}_\epsilon$ solves the elliptic PDE $d^+\triple{a}_\epsilon=\mathcal{F}(d^-\triple{a}_\epsilon\ast d^-\triple{a}_\epsilon)-\triple{\zeta}_\epsilon$, $d^\ast\triple{a}_\epsilon=0$. By elliptic regularity, for any $k\geq 2$ the $C^{k,\alpha}$--norm of $\triple{a}_\epsilon$ on compact sets of $M^\textup{gh}_\epsilon$ and (after rescaling) on compact sets of the gravitational instantons $M_j$ and $N_i$ is controlled in terms of $\| \triple{a}_\epsilon \| _{C^{1,\alpha}_\delta}+\| \triple{\zeta}_\epsilon\|$. In particular, on compact sets of $M^\textup{gh}_\epsilon$ the hyperk\"ahler metric induced by $\triple{\omega}_\epsilon + d\triple{a}_\epsilon + \triple{\zeta}_\epsilon$ is $C^{k,\alpha}$--close to $g_\T + \epsilon^2 \theta^2$. The statements (i), (ii) and (iii) about the limit $\epsilon\ra 0$ now follow.
\endproof
\end{theorem}

By varying all parameters involved in the construction we can in fact realise a whole open set in the moduli space of hyperk\"ahler metrics on the K3 surface. Indeed,
\begin{enumerate}
\item the moduli space of flat tori is $6$--dimensional;
\item the choice of punctures $p_1,\dots, p_n$ yields additional $3n$ parameters;
\item once the punctured torus and weights are fixed, the moduli space of abelian Dirac monopoles with prescribed singularities is $4$ dimensional (one has to choose $\epsilon$ and the $3$--moduli of a flat connection);
\item each $D_{m_j}$ ALF space contributes $3m_j$ parameters and every $A_{k_i-1}$ ALF space contributes $3(k_i-1)$ parameters.
\end{enumerate}
Hence the total number of parameters in the construction is
\[
6+3n+4+3\sum_{j=1}^8{m_j} + 3\sum_{i=1}^n{k_i}-3n=10+3\times 16=58,
\]
which is exactly the dimension of the moduli space of Ricci-flat metrics on the K3 surface (without any normalisation on volume).

\section{Stable minimal surfaces}\label{sec:Minimal:surfaces}

In this final section we exploit our gluing construction of hyperk\"ahler metrics on the K3 surface to deduce some information about their submanifold geometry.

It is well known that holomorphic submanifolds of a K\"ahler manifold are volume minimising. It is a classical problem in the theory of minimal submanifolds in K\"ahler manifolds to understand to what extent volume minimising submanifolds must be holomorphic or anti-holomorphic.

In \cite{Micallef} Micallef showed that every stable minimal surface in a flat $4$--torus must be holomorphic with respect to a complex structure compatible with the metric. (Note however that this is no longer the case for higher dimensional tori \cite{Arezzo:Micallef}.) In view of Micallef's result it was thought for some time that a similar result could hold for the K3 surface endowed with a hyperk\"ahler metric. Partial results in this direction were established by Micallef--Wolfson \cite[Theorem 5.3]{Micallef:Wolfson:1} and motivation for the conjecture came from the fact that, given an arbitrary hyperk\"ahler metric on the K3 surface, every homology class can be represented by the sum of surfaces each of which is holomorphic with respect to some complex structure compatible with the metric. However, Micallef--Wolfson \cite{Micallef:Wolfson:2} have eventually shown that no analogue of the result for $4$--tori holds for the K3 surface. Indeed they found a class $\alpha \in H_2(K3,\Z)$ and a hyperk\"ahler metric $g$ on the K3 surface such that the volume minimiser in $\alpha$ decomposes into a sum of branched minimal surfaces $\Sigma_1 \cup \dots \cup \Sigma_k$ not all of which can be holomorphic with respect to some complex structure compatible with $g$.

We can use our gluing construction to construct further (simpler) examples of strictly stable minimal spheres with respect to some hyperk\"ahler metric on the K3 surface which cannot be holomorphic for any complex structure compatible with the metric.

\begin{theorem}\label{thm:Minimal:Surfaces}
There exist hyperk\"ahler metrics on the K3 surface that contain a strictly stable minimal sphere which is not holomorphic with respect to any complex structure compatible with the metric.
\proof
In \cite[Proposition 5.5]{Micallef:Wolfson:1} Micallef--Wolfson show that the double cover of the Atiyah--Hitchin manifold, the rotationally symmetric $D_1$ ALF space, contains a strictly stable minimal $2$--sphere $\Sigma$ with $[\Sigma]\cdot [\Sigma]=-4$. Since every holomorphic curve $\Sigma$ of genus $\gamma$ in a hyperk\"ahler $4$--manifold must have $[\Sigma]\cdot [\Sigma]=2\gamma -2$ by the adjunction formula, this minimal $2$--sphere cannot be holomorphic with respect to any complex structure. One can also use the isometric action of $SU(2)$ on the Atiyah--Hitchin metric to prove this fact: the $SU(2)$ action preserves the metric but rotates the complex structures (equivalently, the hyperk\"ahler triple) and the minimal $2$--sphere is an $SU(2)$--orbit. Hence the periods $\int_{\Sigma}{\omega_i}$ are forced to vanish.

Now, consider an approximate hyperk\"ahler metric $g_\epsilon$ obtained in Section \ref{sec:Approximate:solution} by using the rotationally symmetric $D_1$ ALF space as one of the building blocks. Thus $g_\epsilon$ contains a strictly stable minimal sphere $\Sigma$ with $[\Sigma]\cdot [\Sigma]=-4$.

Because of strict stability, $\Sigma$ has no Jacobi fields. Then we can invoke White's Implicit Function Theorem for minimal immersions with respect to variations of the ambient metric \cite[Theorem 2.1]{White} to deform $\Sigma$ into a minimal immersion with respect to the hyperk\"ahler metric produced by Theorem \ref{thm:Main:precise} starting from $g_\epsilon$. As before, this minimal $2$--sphere cannot be holomorphic with respect to any complex structure because of its self-intersection number. It is strictly stable by continuity of the spectrum of the Jacobi operator. 
\endproof

The minimal sphere in the double-cover $M$ of the Atiyah--Hitchin manifold has in fact an additional interesting property \cite[Corollary 3.7]{Chen:Li}: it is the image of a harmonic map $u\co S^2 \ra M$ which satisfies \eqref{eqn:Quaternionic:minimal:surfaces}, \ie
\[
du\circ J_{S^2} = -\left( x_1\, J_1\circ du + x_2\, J_2\circ du + x_3\, J_3\circ du \right).
\]
Here $(x_1,x_2,x_3)\co S^2\ra \R^3$ is the standard embedding of $S^2$ as the unit sphere in $\R^3$ and $(J_1,J_2,J_3)$ is a triple of parallel complex structures on $M$ satisfying $J_1 J_2 J_3=-\textup{id}$. Then the radial extension $f$ of $u$ to $\R^3$ is a \emph{Fueter map} $f\co \R^3\ra M$ with an isolated singularity at the origin, \ie it satisfies the Fueter equation
\begin{equation}\label{eqn:Fueter}
J_1 \frac{\partial{f}}{\partial x_1} + J_2 \frac{\partial{f}}{\partial x_2} + J_3 \frac{\partial{f}}{\partial x_3}=0
\end{equation}
on $\R^3\setminus\{ 0\}$. Moreover, the extended map $F\co \R^4 \ra M$ obtained composing $f$ with the projection of $\R^4 = \R^3\times\R$ onto the first factor is \emph{triholomorphic} outside of its $1$--dimensional singular set $\{ 0 \} \times \R$, \ie 
\begin{equation}\label{eqn:Triholomorphic}
dF=J_1 \circ dF \circ I + J_2\circ dF \circ J + J_3\circ dF\circ K.
\end{equation}
Here $I,J,K$ are the standard complex structures on $\R^4\simeq \HH$. Since $f$ and $F$ are invariant under scalings in the domain they provide examples of \emph{tangent maps} for Fueter and triholomorphic maps respectively.

Both \eqref{eqn:Fueter} and \eqref{eqn:Triholomorphic} (and their generalisations to the case of sections of a non-trivial bundle with hyperk\"ahler fibres) arise in different contexts related to gauge theory in higher dimensions, \cf for example \cite{OFarrill} and \cite[\S 6]{Donaldson:Segal}, the definition of enumerative invariants of hyperk\"ahler manifolds \cite{Tian:GW:inv:HK} and the hyperk\"ahler Floer theory of \cite{Salamon:al}. For these geometric applications it is essential to understand the regularity theory of limits of smooth Fueter and triholomorphic maps. In this direction, Bellettini--Tian \cite[Remark 1.2]{Bellettini:Tian} have recently stressed the importance of the existence or non-existence of homogeneous triholomorphic maps with a singular set of codimension $3$. The only known example has a non-compact target, the rotationally symmetric $D_1$ ALF space. The minimal spheres of Theorem \ref{thm:Minimal:Surfaces} give rise to examples of Fueter and triholomorphic tangent maps with a singular set of codimension $3$ and \emph{compact} target.

\begin{theorem}\label{thm:Triholomorphic}
Every non-holomorphic minimal $2$--sphere produced in Theorem \ref{thm:Minimal:Surfaces} is the image of a map $u\co S^2 \ra \textup{K3}$ satisfying \eqref{eqn:Quaternionic:minimal:surfaces}. Moreover the radial extension of $u$ as a Fueter map $f\co \R^3\ra \textup{K3}$ with an isolated singularity at the origin and the extension of $f$ as a homogeneous triholomorphic map $F\co \HH\ra \textup{K3}$ with a singular set of codimension $3$ are both stationary harmonic. 
\end{theorem}
Recall that a stationary harmonic map between two Riemannian manifolds is a critical point for the Dirichlet energy with respect to variations both of the domain and the target.
\proof
The proof of the fact that the minimal $2$--sphere in the double-cover of the Atiayh--Hitchin manifold is the image of a map which satisfies \eqref{eqn:Quaternionic:minimal:surfaces} given in \cite[Corollary 3.7]{Chen:Li} immediately generalises to our case. Indeed, Chen--Li prove that any minimal surface $\Sigma$ in a hyperk\"ahler $4$--manifold $(M^4,J_1,J_2,J_3)$ is the image of a map $u\co \Sigma \ra M$ which satisfies
\[
du\circ J_{\Sigma} = -\left( a_1 J_1 \circ du + a_2 J_2 \circ du + a_3 J_3 \circ du\right)
\] 
for some holomorphic map $\triple{a}=(a_1,a_2,a_3)\co \Sigma \ra S^2\subset\R^3$. Moreover, \cite[Theorem 3.6]{Chen:Li} expresses the degree of $\triple{a}$ in terms of the Euler characteristic and self-intersection number of $\Sigma$ (here we do not distinguish between $\Sigma$ as the source of $u$ and its image in $M$): $2\, \textup{deg}\,\triple{a} +\chi (\Sigma)+[\Sigma]\cdot [\Sigma]=0$. These results immediately imply that up to reparametrization the minimal sphere $\Sigma$ produced by Theorem \ref{thm:Minimal:Surfaces} (for which $[\Sigma]\cdot [\Sigma]=-4$) is the image of a map $u\co S^2\ra \textup{K3}$ satisfying \eqref{eqn:Quaternionic:minimal:surfaces}. In particular the radial extension $f$ of $u$ to $\R^3$ is a Fueter map and the extension of $f$ to $F\co \HH \ra \textup{K3}$ is triholomorphic.

It is well-known (\cf for example \cite[Lemma 3.8]{Chen:Li}) that the extended maps $f$ and $F$ are stationary harmonic if and only if $u$ is balanced, \ie
\[
\int_{S^2}{x_i\, |\nabla u|^2}=0, \qquad \mbox{for } i=1,2,3.
\]
In the case of the double-cover of the Atiyah--Hitchin manifold, Chen--Li prove that the balancing is satisfied by exploiting the invariance under the group $\Z_2$ of deck transformations. We must argue differently. Using \eqref{eqn:Quaternionic:minimal:surfaces} we find
\[
\int_{S^2}{x_i\, |\nabla u|^2}=-2\int_{S^2}{u^\ast\omega_i},
\]  
where $\omega_i$ is the K\"ahler form on the K3 surface corresponding to the complex structure $J_i$. Now, in the proof of Theorem \ref{thm:Minimal:Surfaces} we observed that the minimal sphere in the double-cover of the Atiyah--Hitchin manifold has vanishing periods. Thus any approximate hyperk\"ahler triple $\triple{\omega}_\epsilon$ obtained in Section \ref{sec:Approximate:solution} by using the rotationally symmetric $D_1$ ALF space as one of the building blocks contains a strictly stable minimal sphere $\Sigma$ with $[\Sigma]\cdot[\Sigma]=-4$ and $\int_{\Sigma}{\triple{\omega}_\epsilon}=0$. Since the deformation of $\triple{\omega}_\epsilon=(\omega_1^\epsilon,\omega_2^\epsilon,\omega_3^\epsilon)$ to an exact hyperk\"ahler triple changes the cohomology classes of the three $2$--forms only by constant linear combinations of $[\omega_1^\epsilon], [\omega_2^\epsilon], [\omega_3^\epsilon]$ themselves, we deduce that the periods of the minimal $2$--sphere produced by Theorem \ref{thm:Minimal:Surfaces} all vanish.
\endproof

\end{theorem}
 
\bibliographystyle{amsinitial}
\bibliography{K3}

\end{document}